\documentclass[10pt,a4paper]{article}

\usepackage{amsmath}
\usepackage{graphicx,latexsym,euscript,makeidx,color}
\usepackage{amsmath,amsfonts,amssymb,amsthm,mathrsfs}
\usepackage{enumerate}

\usepackage{amssymb}
\usepackage{amsmath}
\usepackage{amsfonts}
\usepackage{graphicx}
\usepackage{graphicx}          
\usepackage{color}

\usepackage[colorlinks,linkcolor=blue,anchorcolor=green,citecolor=red]{hyperref}

\usepackage{geometry}
\def\sqr#1#2{{\vcenter{\vbox{\hrule height.#2pt
              \hbox{\vrule width.#2pt height#1pt \kern#1pt \vrule width.#2pt}
              \hrule height.#2pt}}}}
\def\signed #1{{\unskip\nobreak\hfil\penalty50
              \hskip2em\hbox{}\nobreak\hfil#1
              \parfillskip=0pt \finalhyphendemerits=0 \par}}
\def\endpf{\signed {$\sqr69$}}

\def\5n{\negthinspace \negthinspace \negthinspace \negthinspace \negthinspace }
\def\4n{\negthinspace \negthinspace \negthinspace \negthinspace }
\def\3n{\negthinspace \negthinspace \negthinspace }
\def\2n{\negthinspace \negthinspace }
\def\1n{\negthinspace }

\def\ms{\medskip}                
               
\def\ds{\displaystyle}

            \def\({\Big (}
                  \def\){\Big )}
          \def\[{\Big[}
           \def\]{\Big]}

\def\bde{\begin{definition}\label}    \def\ede{\end{definition}}
\def\bt{\begin{theorem}\label}        \def\et{\end{theorem}}
\def\bc{\begin{corollary}\label}      \def\ec{\end{corollary}}
\def\bl{\begin{lemma}\label}          \def\el{\end{lemma}}
\def\bp{\begin{proposition}\label}    \def\ep{\end{proposition}}
\def\bas{\begin{assumption}\label}    \def\eas{\end{assumption}}
\def\br{\begin{remark}\label}         \def\er{\end{remark}}
\def\bex{\begin{example}\label}       \def\ex{\end{example}}
\def\ba{\begin{array}}                \def\ea{\end{array}}
\def\be{\begin{equation}}
\def\bel{\begin{equation}\label}      \def\ee{\end{equation}}
\def\bea{\begin{eqnarray*}}           \def\eea{\end{eqnarray*}}

\def\square#1{\vbox{\hrule\hbox{\vrule height#1%
     \kern#1\vrule}\hrule}}
\def\rectangle#1#2{\vbox{\hrule\hbox{\vrule height#1%
     \kern#2\vrule}\hrule}}

\font\tenbb=msbm10 \font\sevenbb=msbm7 \font\fivebb=msbm5
\newfam\bbfam
\scriptscriptfont\bbfam=\fivebb \textfont\bbfam=\tenbb
\scriptfont\bbfam=\sevenbb

\newtheorem{theorem}{\indent Theorem}[section]
\newtheorem{definition}[theorem]{\indent Definition}
\newtheorem{proposition}[theorem]{\indent Proposition}
\newtheorem{corollary}[theorem]{\indent Corollary}
\newtheorem{lemma}[theorem]{\indent Lemma}
\newtheorem{remark}[theorem]{\indent Remark}
\newtheorem{example}[theorem]{\indent Example}

\newtheorem{assumption}[theorem]{\indent Assumption}

\makeatletter
   
   \@addtoreset{equation}{section}
\makeatother

\allowdisplaybreaks[4]
\def\hpm{\hphantom}

\begin{document}

\title{\bf Delayed Optimal Control of Stochastic LQ Problem
\thanks{Part of this paper was presented at the 36th Chinese Control Conference.}
}
\author{Yuan-Hua Ni\thanks{College of Computer and Control Engineering,  Nankai University, Tianjin 300350 , P. R. China}~~~~ Cedric Ka-Fai Yiu\thanks{Department of
Applied Mathematics, The Hong Kong Polytechnic University, Hunghom,
Kowloon, Hong Kong, P.R. China.}~~~~ Huanshui Zhang\thanks{School of Control Science and Engineering, Shandong University, Jinan 250061, P.R. China.}~~~~ Ji-Feng Zhang\thanks{Key Laboratory of Systems and
Control, Institute of Systems Science, Academy of Mathematics and
Systems Science, Chinese Academy of Sciences, and the School of Mathematical Sciences, University of Chinese Academy of Sciences,
P. R. China.}}
\maketitle

{\bf Abstract:} A stochastic LQ problem with multiplicative noises and transmission delay is studied in this paper, which does not require any definiteness constraint on the cost weighting matrices.  From some abstract representations of the system and cost functional, the solvability of this LQ problem is characterized by some conditions with operator form. Based on these, necessary and sufficient conditions are derived for the case with a fixed time-state initial pair and the general case with all the time-state initial pairs.
For both cases, a set of coupled discrete-time Riccati-like equations
can be derived to characterize the existence and the form of the delayed optimal control. In particular, for the general case with all the initial pairs, the existence of the delayed optimal control is equivalent to the solvability of the Riccati-like equations with some algebraic constraints, and both of them are also equivalent to the solvability of a set of coupled linear matrix equality-inequalities. Note that both the constrained Riccati-like equations and the linear matrix equality-inequalities are introduced for the first time in the literature for the proposed LQ problem. Furthermore, the convexity and the {uniform convexity} of the cost functional are fully characterized via certain properties of the solution of the Riccati-like equations.

\ms


\ms


\section{Introduction}

Linear-quadratic (LQ, for short) optimal control was pioneered by Kalman \cite{Kalman} in 1960, which is now a classical yet fundamental problem in control theory.
%
%
Extension to stochastic LQ problems was first carried out by Wonham \cite{Wonham} in 1968, and has received considerable interests and efforts since then. 
A common assumption of most literature on stochastic LQ problems is that the state weighting matrices are nonnegative definite and the control weighting matrices are positive definite. 
Contrary to this, Chen, Li and Zhou \cite{Chen-Li-Zhou} revealed in 1998 that a stochastic LQ problem with multiplicative noises might still be solvable even if the cost weighting matrices are indefinite. More about this kind of LQ problems can be found in \cite{Ait-Chen-Zhou-2002} \cite{Ait-Moore-Zhou} \cite{Hu-Zhou-2003} \cite{Sun-Yong-Siam-2016} and references therein. Recently, some researchers are  interested in the so-called mean-field LQ problems \cite{Ni-Elliott-Li} \cite{Ni-Zhang-Li-2015} \cite{Sun-Yong-2016} \cite{Wang-guanchen} \cite{Yong-2013} \cite{Yong-2015incon}.
An important feature of mean-field control problems is that the expected values of the state and control enter nonlinearly into the cost functional, which will bring new phenomena and new theoretical difficulties.

Note that all the aforementioned papers are free of time delay. If time delay happens to appear in the system state, the control input or the information-transmission channel, it is much more complicated and challenging to design the optimal control of the corresponding LQ problems. Such kind of LQ problems have been extensively studied since 1970's; see, for example, \cite{Alekal-Lee1971} \cite{Delfour1986} \cite{Koivo-Lee1972} \cite{Watanabe-Ito-1981} \cite{Zhang-huanshui2015} or other related literature \cite{Altman-Basar-Srikant-1999} \cite{Kojima-Ishijima-2006} \cite{Mao-Deng-Wan-2016} \cite{Tadmor-Mirkin-2005}. Concerned with a deterministic LQ problem with input delay, it is shown \cite{Watanabe-Ito-1981} that the delayed optimal control is obtained by invoking the Smith predictor theory, and that the optimal gains are same to those of the LQ problem without input delay.
Unfortunately, the results about deterministic LQ problems (with input delay) cannot be directly generalized to the stochastic setting. In \cite{Zhang-huanshui2015}, the authors considered a discrete-time stochastic LQ problem with input delay and multiplicative noises, and showed that the optimal control (if exists) is a linear feedback of $d$-step-lagged conditional expectation of current states and that the optimal gains are computed via a set of coupled discrete-time Riccati-like equations. Here, the set of discrete-time Riccati-like equations differs significantly from what we have in hand the standard discrete-time Riccati equation.

It is worth pointing out that the stochastic systems with multiplicative noises have been extensively studied in the past half century.
From the viewpoint of mathematics, almost all the theories about stochastic differential equations (SDEs, for short) are for the case with multiplicative noises, and there are lots of practical motivations to study such kind of SDEs.
The study of controlled systems with multiplicative noises is also popular in the control community; a recent small collection in the literature related to our paper includes \cite{Ait-Chen-Zhou-2002} \cite{Ait-Moore-Zhou} \cite{Bismut} \cite{Chen-Wu-2010} \cite{Chen-Li-Zhou} \cite{Hu-Zhou-2003} \cite{Huang-jianhui} \cite{Ni-Elliott-Li} \cite{Peng} \cite{Sun-Yong-Siam-2016} \cite{Wang-Yang-Yong-Yu} \cite{Yong-2013}.

In this paper, a general discrete-time stochastic LQ problem with multiplicative noises and transmission delay is thoroughly investigated, whose cost weighting matrices for the state and control are allowed to be indefinite. Apart from intending to generalize the existing results \cite{Alekal-Lee1971} \cite{Delfour1986} \cite{Koivo-Lee1972} \cite{Watanabe-Ito-1981} \cite{Zhang-huanshui2015} to the joint case with indefiniteness and time delay, the topic of this paper is also partially motivated by recent progresses in network control system and other related areas. 
Transmission delay, or sometimes called as communication delay, is a key feature of network control systems \cite{Baillieul-Antsaklis} \cite{Hespanha} \cite{Wei-Zhang-Fu}, which is generally caused by the limited bit rate of communication channels. In fact, transmission delay has been extensively studied in the areas such as discrete-event dynamic systems \cite{Delay-Wonham}, multi-agent systems \cite{Delay-Xie} \cite{Delay-Tian} \cite{Delay-Ren}, networked mobile robots \cite{Delay-Robot}, receding horizon control \cite{Delay-Receding}, flexible spacecraft \cite{Delay-Du}, and so on. Furthermore, such kind of delays are also related to the measurement delays \cite{Dealy-Cacace} \cite{Delay-Cao} \cite{Delayed-Nagpal} \cite{Delay-meausrenent-4} \cite{Delay-measurement-zhang}, which arise in measurement channels.

The contents of this paper are as follows. For the completeness and parallel to that in \cite{Yong-Zhou}, the considered problem (Problem (LQ)) is converted in Section 3 to a quadratic optimization problem in the Hilbert space. By this reformulation, we can derive some abstract conditions on the solvability of Problem (LQ), which gives us an overall perspective of Problem (LQ) and motivates the analysis of the sections followed. This part of work is a discrete-time version (with   state transmission delay) of the results in \cite{Yong-Zhou}, and the backward stochastic difference equations (BS$\Delta$Es, for short) are involved here.

In Section 4, for the case with a fixed time-state initial pair, the solvability of Problem (LQ) at that initial pair is equivalent to that a stationary condition and a convexity condition are satisfied, with the backward state of a forward-backward stochastic difference equation (FBS$\Delta$E, for short) being involved in the stationary condition. Further, a set of coupled discrete-time Riccati-like equations is introduced, by which we can express the backward state of the FBS$\Delta$E via its forward state.
     Moreover, equivalent characterizations of the stationary condition and the convexity condition are derived via certain properties of the solution of the Riccati-like equations.

In Section 5, for the case with all the time-state  initial pairs, the following facts are shown to be equivalent: (i) Problem (LQ) is finite; (ii) Problem (LQ) is solvable; (iii) a set of constrained coupled discrete-time Riccati-like equations is solvable; (iv) a set of coupled linear matrix equality-inequalities (LMEIs, for short) is solvable.
Moreover, the unique solvability of Problem (LQ) at the initial pair $(t,x)$ is shown to be  equivalent to the unique solvability at \emph{any} initial pair $(k,\xi)\in \{t,...,N-1\}\times \mathbb{R}^n$, both of which are equivalent to the uniform convexity of the cost functional and the positive definiteness of certain matrices involved in the constrained Riccati-like equations.

From our derived results, we have the following remarks.
\begin{itemize}

\item For Problem (LQ), the case with a fixed time-state initial pair differs significantly from the case with all the time-state  initial pairs; this can be seen from Theorem \ref{Theorem--Nece-Suff-fixed-final} and Theorem \ref{Theorem-all-ii}. Hence, we separately discuss the two cases.

\item By the stationary condition and a backward procedure of calculations, we can get the Riccati-like equations (\ref{Riccati-1})-(\ref{Riccati-3}) and express FBS$\Delta$E's backward state via its forward state and the solution of the Riccati-like equations. Due to the $d$-step-lagged information structure, the Riccati-like equations are much more complicated than the standard discrete-time Riccati equation.

\item The convexity of the cost functional is fully characterized in Theorem \ref{Theorem-convex-nece-suff} via certain properties of solution of the Riccati-like equations (\ref{Riccati-1})-(\ref{Riccati-3}), which is proved by using a technique of control shifting.
To the best of our knowledge, this result seems to be the first one of equivalent characterization on the convexity of the cost functional of LQ problem.

Based on this, necessary and sufficient conditions on the solvability of Problem (LQ) for a fixed initial pair is presented.

\item Note that the constrained Riccati-like equations (\ref{Riccati-1-2})-(\ref{Riccati-3-2}) and the LMEIs (\ref{Riccati-inequality-1})-(\ref{Riccati-inequality-3}) are introduced for the first time, to the best of our knowledge. Furthermore, from a solution of the LMEIs, an explicit procedure is presented to construct a solution of the constrained Riccati-like equations. Such a procedure is potentially useful to study the algebraic Riccati-like equations that we will encounter in the infinite-horizon version of Problem (LQ).

It is worth mentioning that there are linear equations in the set of Riccati-like equations and the LMEIs contain equality constraints. Note that such new feature do not appear in deterministic LQ problems (with time delay) and standard stochastic LQ problems.

\end{itemize}

%

%

In \cite{Zhang-huanshui2015}, stochastic LQ problems with multiplicative noises and input delay were investigated, whose  cost weighting matrices are assumed to be nonnegative definite.
%
%
%
This paper is of general indefinite case, and thus, differs substaintially from \cite{Zhang-huanshui2015}. In the context of this paper,
it is proved in \cite{Zhang-huanshui2015} that (ii) and (v) of Theorem \ref{Theorem-nec-suff} are equivalent for the nonnegative-definite case, which is
the main result of the finite-horizon LQ problem in \cite{Zhang-huanshui2015}.
%
%
%
Note that in \cite{Zhang-huanshui2015}, the case with a fixed initial pair and the case with all the initial pairs are not differentiated, and no LMEIs are mentioned.
Hence, the results of this paper are broader than those of the finite-horizon LQ problem of \cite{Zhang-huanshui2015}.
Furthermore, the transmission delay is studied in this paper, which is different from the input delay \cite{Zhang-huanshui2015}; this is why the Riccati-like equations of this paper are divided into several pieces.

The rest of this paper is organized as follows. Section
\ref{Section-2} and Section \ref{Section-3} give the problem formulation and an abstract consideration. In Section \ref{section--fixed} and Section \ref{section--all}, the case with a fixed initial pair and the case with all the initial pairs are investigated, respectively.
Section \ref{Section-example} gives an example, and some concluding remarks are given in Section \ref{section-conclusion}.

\section{Problem formulation}\label{Section-2}

Consider the following controlled stochastic difference equation (S$\Delta$E, for short)
\begin{eqnarray}\label{system-1}
\left\{\begin{array}{l}
X_{k+1}=\big{(}A_{k}X_k+B_{k}u_{k}\big{)}+\big{(}C_{k}X_k+D_{k}u_{k}\big{)}w_k, \\[1mm]
X_t=x,~~k\in \mathbb{T}_t\triangleq \{t,t+1,...,N-1\},~~t\in \mathbb{T}\triangleq \{0,1,...,N-1\},
\end{array}
\right.
\end{eqnarray}
where $A_{k}, C_{k}\in \mathbb{R}^{n\times
n}$, $B_{k}, D_{k}\in\mathbb{R}^{n\times m}$ are deterministic matrices. The noise $\{w_k, k\in
\mathbb{T}\}$ is assumed to be  a martingale difference sequence
defined on a probability space $(\Omega, \mathcal{F}, P)$ with
\begin{eqnarray}\label{w-moment}
\mathbb{E}_{k+1}[w_{k+1}]=0,~~\mathbb{E}_{k+1}[(w_{k+1})^2]=1,~k\in \mathbb{T}.
\end{eqnarray}
Here, $\mathbb{E}_{k+1}$ is the conditional mathematical expectation
$\mathbb{E}[\,\cdot\,|\mathcal{F}_{k+1}]$ with respect to $\mathcal{F}_{k+1}=\sigma\{w_l, l=0, 1,\cdots,k\}$, and $\mathcal{F}_{0}$ is understood as
$\{\emptyset, \Omega\}$.
Introduce the following cost functional associated with (\ref{system-1})
\begin{eqnarray}\label{cost-1}
&&\hspace{-2em}J(t,x; u)
 =\sum_{k=t}^{N-1}\mathbb{E}\big{[}X_k^TQ_{k}X_k+u_{k}^TR_{k}u_{k}\big{]}+\mathbb{E}\big{[}X_N^TGX_N\big{]},
\end{eqnarray}
where $Q_{k}, R_{k}, k\in \mathbb{T}_t$, $G$ are deterministic symmetric matrices
of appropriate dimensions. Note, here, that we do not pose any definiteness constraints on the cost weighting matrices.

 %
%
This paper is concerned with the case with transmission delay. For such kinds of time delays and the related measurement delays, find  \cite{Dealy-Cacace} \cite{Delay-Cao}  \cite{Baillieul-Antsaklis} \cite{Delay-Du} \cite{Hespanha} \cite{Delay-Receding} \cite{Delay-Ren} \cite{Delay-Xie} \cite{Delayed-Nagpal}  \cite{Delay-meausrenent-4} \cite{Delay-Tian} \cite{Wei-Zhang-Fu} \cite{Delay-Robot} \cite{Delay-Wonham}  \cite{Delay-measurement-zhang} in Introduction for their motivations and applications.
Assume in this paper that there is a $d$-step time delay in the transmission/measurement channel ($d\geq 2$).
%
%
%
Due to this, 
for $k\in \{t,...t+d\}$ no new information is available and the controller's decision information set remains $\mathcal{F}_t$; and for $k\in \mathbb{T}_{t+d}=\{t+d,...,N-1\}$ the information set should be $\mathcal{F}_{k-d}$. 
In this paper, we select 
\begin{eqnarray}\label{admissible-control--t}
\mathcal{U}_{ad}^t=\big{(}l^2_\mathcal{F}(t; \mathbb{R}^m)\big{)}^{d}\times l^2_\mathcal{F}(\mathbb{T}_{t}^{-d}; \mathbb{R}^m)
\end{eqnarray}
as the admissible control set, where 
\begin{eqnarray}\label{L^2-t}
l^2_\mathcal{F}(t; \mathbb{R}^m)=\Big{\{}\zeta \in \mathbb{R}^m\,\big{|}\, \zeta\mbox{ is }\mathcal{F}_t\mbox{-measurable, and }\mathbb{E}|\zeta|^2<\infty \Big{\}},~~t=0,...,N,
\end{eqnarray}
and
\begin{eqnarray}\label{L^2-t-d}
l^2_\mathcal{F}(\mathbb{T}^{-d}_{t}; \mathbb{R}^m)=\Big{\{}\nu=\{\nu_k, k\in \mathbb{T}_{t}^{-d}\}\,\big{|}\,\nu_k\mbox{ is }\mathcal{F}_k\mbox{-measurable, and }\mathbb{E}|\nu_k|^2<\infty, k\in\mathbb{T}_{t}^{-d} \Big{\}}
\end{eqnarray}
with $$\mathbb{T}_{t}^{-d}=\{t,...,N-1-d\}.$$
Therefore, for any $(u_t,...,u_{N-1})=u\in \mathcal{U}_{ad}^t$, $u_k$ is $\mathcal{F}_t$-measurable if $k\in \{t,...t+d\}$, and $u_k$ is $\mathcal{F}_{k-d}$-measurable if $k\in \mathbb{T}_t^{-d}$; this reflects the property of causality.

The following optimal control problem will be studied in this paper.

\textbf{Problem (LQ).} \emph{For a time-state initial pair $(t,x)\in \mathbb{T}\times l^2_\mathcal{F}(t; \mathbb{R}^n) $, find a
$\bar{u}\in \mathcal{U}_{ad}^t$ such that
\begin{eqnarray}\label{Problem-LQ}
J(t,x;\bar{u}) = \inf_{u\in \mathcal{U}_{ad}^t}J(t,x;u).
\end{eqnarray}}

\begin{remark}
Noting that the initial pair $(t,x)$ is specialized, hereafter the above problem will be called as Problem (LQ) for the initial pair $(t,x)$.
Furthermore, any $\bar{u}$ satisfying (\ref{Problem-LQ}) is called an optimal control  of Problem (LQ) for the initial pair $(t,x)$.

\end{remark}

%

\begin{definition}
Problem (LQ) is said to be (uniquely) solvable at $(t,x)$ if there exists a (unique) $\bar{u}\in \mathcal{U}^t_{ad}$
such that (\ref{Problem-LQ}) holds. 
\end{definition}

In what follows, we shall review some knowledge on matrix.
Recall the pseudo-inverse of a matrix. By \cite{Penrose}, for a
given matrix $M\in \mathbb{R}^{n\times m}$, there exists a unique
matrix in $\mathbb{R}^{m\times n}$ denoted by $M^\dagger$ such that
\begin{eqnarray}
\left\{
\begin{array}{l}
MM^\dagger M=M,~~ M^\dagger M M^\dagger=M^\dagger,\\
(MM^\dagger)^T=MM^\dagger, ~~(M^\dagger M)^T=M^\dagger M.
\end{array}
\right.
\end{eqnarray}
This $M^\dagger$ is called the Moore-Penrose inverse  of $M$. The
following lemma is from \cite{Ait-Chen-Zhou-2002}.

\begin{lemma}\label{Lemma-matrix-equation}
Let matrices $L$, $M$ and $N$ be given with appropriate size. Then,
$LXM=N$ has a solution $X$ if and only if $LL^\dagger NMM^\dagger=N$.
Moreover, the solution of $LXM=N$ can be expressed as
$X=L^\dagger NM^\dagger+Y-L^\dagger LYMM^\dagger$,
where $Y$ is a matrix with appropriate size.
\end{lemma}

If $M=I$ in Lemma \ref{Lemma-matrix-equation}, then $LL^\dagger N=N$ is equivalent to $\mbox{Ran}(N)\subset \mbox{Ran}(L)$. Here, $\mbox{Ran}(N)$ is the range of $N$. The following is the so-called extended Schur's lemma.

\begin{lemma}\label{Lemma-Schur}
Let $S=S^T\in \mathbb{R}^{n\times n}, W=W^T\in \mathbb{R}^{m\times m}, H\in \mathbb{R}^{m\times n}$. Then
\begin{eqnarray*}
\left[
\begin{array}{cc}
S&H^T\\H&W
\end{array}
\right]\geq 0
\end{eqnarray*}
if and only if
\begin{eqnarray*}
S-H^TW^\dagger H\geq 0,~W\geq 0,~WW^\dagger H=H.
\end{eqnarray*}

\end{lemma}

\section{An abstract consideration}\label{Section-3}

For the completeness of theory, in this section, we convert Problem (LQ) to a quadratic optimization problem in Hilbert space, based on which some necessary conditions and sufficient conditions are given on the solvability of Problem (LQ).
This part of work is a discrete-time version (with   state transmission delay) of the results in \cite{Yong-Zhou}, which will give us an overall perspective of Problem (LQ) and will motivate the analysis of the following sections.

To begin with, for $k,l \in \mathbb{T}_t$, let
\begin{eqnarray*}
\left\{\begin{array}{l}
{\Phi}(k,\ell)=(A_{k}+w_{k}C_{k})(A_{k-1}+w_{k-1}C_{k-1})\cdots (A_{\ell}+w_{\ell}C_{\ell}),~~ k>\ell, \\[1mm]
\Phi(k,k)=A_{k}+w_{k}C_{k},\\[1mm]
{\Phi}(k,\ell)=I, k<\ell.
\end{array}
\right.\end{eqnarray*}
From (\ref{system-1}), we have
\begin{eqnarray}\label{system-operator-0}
X_{k+1} ={\Phi}(k,t)x +\ds\ds\sum_{\ell=t}^{k}\Phi(k,\ell+1)(B_{\ell}+w_{\ell}D_{\ell})u_{\ell},~~k\in \mathbb{T}_t.
\end{eqnarray}
For any $x\in l^2_{\mathcal{F}}(t; \mathbb{R}^n)$ and $u\in \mathcal{U}^t_{ad}$, define the following operators
\begin{eqnarray*}
\left\{\begin{array}{l}
\Gamma^tx=\Big{\{}\big{(}(\Gamma^tx)_t,...,(\Gamma^tx)_{N-1}\big{)}\,\Big{|}\,(\Gamma^t x)_k={\Phi}(k-1,t)x,~k\in \mathbb{T}_t\Big{\}}, \\ [2mm]
\hat{\Gamma}^t x={\Phi}(N-1,t)x, \\
L^tu=\Big{\{}\big{(}(L^tu)_t,...,(L^tu)_{N-1}\big{)}\,\Big{|}\,(L^t u)_t=0, (L^t u)_k=\ds\sum_{\ell=t}^{k-1}\Phi(k-1,\ell+1)(B_{\ell}+w_{\ell}D_{\ell})u_{\ell},~k\in \mathbb{T}_{t+1}\Big{\}},\\
\hat{L}^tu=\ds\sum_{\ell=t}^{N-1}\Phi(N-1,\ell+1)(B_{\ell}+w_{\ell}D_{\ell})u_{\ell}.
\end{array}\right.
\end{eqnarray*}
Hence,
\begin{eqnarray}\label{system-operator}
X_k=(\Gamma^t x)_k+(L^tu)_k, ~~k\in {\mathbb{T}}_{t},
\end{eqnarray}
and
\begin{eqnarray}\label{system-operator-2}
X_{N}=\hat{\Gamma}^t x+\hat{L}^tu.
\end{eqnarray}
It is not hard to see that the operators
\begin{eqnarray}\label{operators}
\left\{\begin{array}{l}
\Gamma^t: l^2_{\mathcal{F}}(t; \mathbb{R}^n)\mapsto l^2_{\mathcal{F}}({\mathbb{T}}_{t}; \mathbb{R}^n),~~\hat{\Gamma}^t: l^2_{\mathcal{F}}(t; \mathbb{R}^n)\mapsto l^2_{\mathcal{F}}(N; \mathbb{R}^n),\\
L^t: \mathcal{U}^t_{ad}\mapsto l^2_{\mathcal{F}}({\mathbb{T}}_{t}; \mathbb{R}^n),~~\hat{{L}}^t: \mathcal{U}^t_{ad}\mapsto l^2_{\mathcal{F}}(N; \mathbb{R}^n)
\end{array}\right.
\end{eqnarray}
are all bounded and linear.
Notice that the spaces in (\ref{operators}) are all Hilbert spaces. Therefore, the corresponding adjoint operators uniquely exist.
For $\eta\in l^2_{\mathcal{F}}(N; \mathbb{R}^n)$ and $\xi\in l^2_{\mathcal{F}}({{\mathbb{T}}}_{t}; \mathbb{R}^n)$, introduce the following BS$\Delta$E
\begin{eqnarray}\label{adjoint-1}
\left\{
\begin{array}{l}
%
V_k=A_k^T\mathbb{E}_kV_{k+1}+C_k^T\mathbb{E}_k(V_{k+1}w_k)+\xi_k,\\[1mm]
V_N=\eta,~~k\in \mathbb{T}_t.
\end{array}
\right.
\end{eqnarray}

\begin{proposition}
Let $V^0$ be the solution of (\ref{adjoint-1}) with $\eta=0$ and $V^{00}$ be the solution of (\ref{adjoint-1}) with $\xi=0$. Then the adjoint operators $\Gamma^{t*}$, $L^{t*}$, $\hat{\Gamma}^{t*}$ and $\hat{L}^{t*}$
are given, respectively, by
\begin{eqnarray}\label{adjoint-G}
&&\Gamma^{t*}\xi=V_t^0,\\[1.5mm]
&&\label{adjoint-L}
(L^{t*}\xi)_k=
%
B_k^T\mathbb{E}_{k-d}V^0_{k+1}+D_k^T\mathbb{E}_{k-d}(V^0_{k+1}w_k),~~k\in \mathbb{T}_{t},
\\[1.5mm]
&&\label{adjoint-hat-G}
\hat{\Gamma}^{t*}\eta= V^{00}_t,
\end{eqnarray}
and
\begin{eqnarray}\label{adjoint-hat-L}
(\hat{L}^{t*}\eta)_k=
%
B_k^T\mathbb{E}_{k-d}V^{00}_{k+1}+D_k^T\mathbb{E}_{k-d}(V^{00}_{k+1}w_k),~~k\in \mathbb{T}_{t}.
\end{eqnarray}
In (\ref{adjoint-L}), (\ref{adjoint-hat-L}), $\mathbb{E}_{k-d}$ is understood as $\mathbb{E}_t$ if $k\in \{t,...,t+d-1\}$ (i.e., $k-d<t$).

\end{proposition}

\emph{Proof}. From (\ref{adjoint-1}) and by substituting $X_{k+1}$, we have
\begin{eqnarray}\label{adjoint-2}
&&\hspace{-2em}\mathbb{E}\big{[}\eta^TX_N-V_t^Tx\big{]}=\sum_{k=t}^{N-1}\mathbb{E}\Big{[}V_{k+1}^TX_{k+1}-V_{k}^TX_{k}  \Big{]}\nonumber  \\
&&\hspace{-2em}\hphantom{\mathbb{E}\big{[}\eta^TX_N-V_t^Tx\big{]}}=\sum_{k=t}^{N-1}\mathbb{E}\Big{[}\big{(}A_k^TV_{k+1}+C_k^T(V_{k+1}w_k)-V_k \big{)}^TX_k \Big{]}  \nonumber \\
&&\hspace{-2em}\hphantom{\mathbb{E}\big{[}\eta^TX_N-V_t^Tx\big{]}=}+\sum_{k=t}^{N-1}\mathbb{E}\Big{[}\big{(}B_k^TV_{k+1}+D_k^T(V_{k+1}w_k) \big{)}^Tu_{k} \Big{]}\nonumber\\
&&\hspace{-2em}\hphantom{\mathbb{E}\big{[}\eta^TX_N-V_t^Tx\big{]}}=-\sum_{k=t}^{N-1}\mathbb{E}\big{[}\xi_k^TX_k\big{]}+\sum_{k=t}^{N-1}\mathbb{E}\Big{[}\big{(}B_k^T\mathbb{E}_{k-d}V_{k+1}+D_k^T\mathbb{E}_{k-d}(V_{k+1}w_k) \big{)}^Tu_{k} \Big{]}.
%
\end{eqnarray}
%
%
%
%
Letting $\eta=0, u=0$ in (\ref{adjoint-2}), from (\ref{system-operator}) we have
\begin{eqnarray*}
\langle \Gamma^t x, \xi\rangle=\sum_{k=t}^{N-1} \mathbb{E}\big{[}\xi_k^T(\Gamma^t x)_k\big{]}=\sum_{k=t}^{N-1} \mathbb{E}\big{[}\xi_k^TX_k\big{]}=\mathbb{E}\big{[}x^TV_t^0 \big{]}=\langle x, V^0_t\rangle,
\end{eqnarray*}
which implies (\ref{adjoint-G}). Letting $x=0, \eta=0$ in (\ref{adjoint-2}), the following equation holds
\begin{eqnarray*}
\langle L^tu, \xi \rangle=\sum_{k=t}^{N-1} \mathbb{E}\big{[}(L^tu)_k^T\xi_k \big{]}=\sum_{k=t}^{N-1} \mathbb{E}\big{[}X_k^T\xi_k \big{]}=\sum_{k=t}^{N-1}\mathbb{E}\Big{[}u_{k}^T\big{(}B_k^T\mathbb{E}_{k-d}V^0_{k+1}+D_k^T\mathbb{E}_{k-d}(V^0_{k+1}w_k) \big{)} \Big{]}.
\end{eqnarray*}
Hence, the adjoint operator $L^{t*}$ of $L$ is given by (\ref{adjoint-L}).

Letting $\xi=0, u=0$ in (\ref{adjoint-2}), we have
\begin{eqnarray*}
\langle\hat{\Gamma}^t x, \eta\rangle= \mathbb{E}\big{[}\eta^TX_N \big{]}=\mathbb{E}\big{[}x^TV^{00}_t \big{]}=\langle x, V^{00}_t\rangle.
\end{eqnarray*}
Then, the adjoint operator $\hat{\Gamma}^{t*}$ of $\hat{\Gamma}^t$ is given by (\ref{adjoint-hat-G}). Furthermore, letting $\xi=0, x=0$ in  (\ref{adjoint-2}), it holds that
\begin{eqnarray*}
\langle \hat{L}^tu, \eta\rangle= \mathbb{E}\big{[}\eta^TX_N \big{]}=\sum_{k=t}^{N-1}\mathbb{E}\Big{[}u_{k}^T\big{(}B_k^T\mathbb{E}_{k-d}V^{00}_{k+1}+D_k^T\mathbb{E}_{k-d}(V^{00}_{k+1}w_k) \big{)} \Big{]}.
\end{eqnarray*}
We therefore have (\ref{adjoint-hat-L}). \endpf 

We further use the convention
\begin{eqnarray*}
\left\{\begin{array}{l}
(QX)_k=Q_kX_k,~ k\in \mathbb{T}_t,~\forall X\in l^2_{\mathcal{F}}({T}_t; \mathbb{R}^n),\\
(Ru)_k=R_ku_k,~ k\in \mathbb{T}_t,~\forall u\in \mathcal{U}^t_{ad}.
\end{array}\right.
\end{eqnarray*}
Then, the cost functional $J(t,x; u)$ has the following form
\begin{eqnarray}\label{operator-J}
&&J(t,x; u)=  \langle Q(\Gamma^t x+L^tu), \Gamma^t x+L^tu\rangle_{l^2_{\mathcal{F}}({T}_t; \mathbb{R}^n)}+\langle Ru, u\rangle_{\mathcal{U}^t_{ad}}+\langle G(\hat{\Gamma}^t x+\hat{L}^tu), \hat{\Gamma}^t x+\hat{L}^tu\rangle_{l^2_{\mathcal{F}}(N; \mathbb{R}^n)}\nonumber \\
&&\hphantom{J(t,x; u)}= \langle \Theta^t_1u, u\rangle_{\mathcal{U}^t_{ad}}+2\langle \Theta^t_2x, u\rangle_{\mathcal{U}^t_{ad}}+\langle \Theta^t_3 x, x\rangle_{l^2_{\mathcal{F}}(t; \mathbb{R}^n)}
\end{eqnarray}
with
\begin{eqnarray*}
\left\{
\begin{array}{l}
\Theta^t_1=R+L^{t*}QL^t+\hat{L}^{t*}G\hat{L}^t, \\ [1mm]
\Theta^t_2=L^{t*}Q\Gamma^t+\hat{L}^{t*}G\hat{\Gamma}^t, \\ [1mm]
\Theta^t_3=\Gamma^{t*}Q\Gamma^t+\hat{\Gamma}^{t*}G\hat{\Gamma}^t.
\end{array}
\right.
\end{eqnarray*}
In (\ref{operator-J}), the inner products are for different Hilbert spaces.
%
%
%
%

Based on above preparations, we have the following result.

\begin{proposition}\label{proposition-operator-Hilbert}
The following statements hold.


(i) Let $u, v\in \mathcal{U}_{ad}^t$ and $\lambda\in \mathbb{R}$. Then
\begin{eqnarray*}\label{}
J(t,x;u+\lambda v)-J(t,x;u)=\lambda^2 \langle  \Theta_1^tv, v \rangle_{\mathcal{U}^t_{ad}}+2\lambda \langle \Theta_1^tu+\Theta_2^tx, v\rangle_{\mathcal{U}^t_{ad}}.
\end{eqnarray*}

(ii) Problem (LQ) is (uniquely) solvable at $(t,x)$ if and only if $\Theta^t_1\geq 0$ and there exists a (unique) $\bar{u}\in \mathcal{U}_{ad}^t$ such that
\begin{eqnarray*}
\Theta_1^t\bar{u}+\Theta_2^t x=0.
\end{eqnarray*}

(iii) If $\Theta^t_1>aI$ for some $a>0$, then $J(t,x; u)$ admits a unique minimizer $\bar{u}$
\begin{eqnarray*}\label{proposition-operator-Hilbert-control}
\bar{u}_k=-((\Theta_1^t)^{-1}\Theta^t_2x)_k, ~~k\in\mathbb{T}_t.
\end{eqnarray*}
In addition, if
\begin{eqnarray}\label{proposition-operator-Hilbert-condition}
Q_k\geq 0,~~R_k>0,~k\in\mathbb{T}_t,~G\geq 0,
\end{eqnarray}
then, $\Theta^t_1>aI$ for some $a>0$.
\end{proposition}

\emph{Proof.} (i) follows from (\ref{operator-J}), which implies (ii) and (iii).\endpf
%

%
%
Some calculations show
\begin{eqnarray}\label{operator-theta-1}
(\Theta^t_1u)_k=R_ku_k+B_k^T\mathbb{E}_{k-d}V^1_{k+1}+D_k^T\mathbb{E}_{k-d}(V^1_{k+1}w_k{)},~k\in \mathbb{T}_{t},
\end{eqnarray}
and
\begin{eqnarray}\label{operator-theta-2}
(\Theta^t_2x)_k=
B_k^T\mathbb{E}_{k-d}V^2_{k+1}+D_k^T\mathbb{E}_{k-d}(V^2_{k+1}w_k{)},~k\in \mathbb{T}_{t},
\end{eqnarray}
where $V^1, V^2$ are given by
\begin{eqnarray}\label{adjoint-V-1}
\left\{
\begin{array}{l}
V^1_k=A_k^T\mathbb{E}_kV^1_{k+1}+C_k^T\mathbb{E}_k(V^1_{k+1}w_k)+(QL^tu)_k,\\[1mm]
V^1_N=G\hat{L}^tu,~~k\in \mathbb{T}_t,
\end{array}
\right.
\end{eqnarray}
and
\begin{eqnarray}\label{adjoint-V-2}
\left\{
\begin{array}{l}
V^2_k=A_k^T\mathbb{E}_kV^2_{k+1}+C_k^T\mathbb{E}_k(V^2_{k+1}w_k)+(Q\Gamma^tx)_k,\\[1mm]
V^2_N=G\hat{\Gamma}^tx,~~k\in \mathbb{T}_t.
\end{array}
\right.
\end{eqnarray}
Hence, we have the following results.

\begin{corollary}\label{Lemma-difference}
Let $u, {v}\in \mathcal{U}_{ad}^t$ and  $\lambda \in
\mathbb{R}$. Then,
\begin{eqnarray}\label{Lemma-difference-0}
{J}(t,x; u+\lambda {v})-{J}(t,x; u)=\lambda^2 J(t,0; {v})+2\lambda \sum_{k=t}^{N-1}\mathbb{E}\Big{[}\big{(}R_k u_k+B_k^TZ_{k+1}+D_k^TZ_{k+1}w_k\big{)}^T{v}_k \Big{]},
\end{eqnarray}
where
\begin{eqnarray*}\label{Lemma-difference-Z}
\left\{
\begin{array}{l}
Z_k=Q_k X_k+A_k^T\mathbb{E}_k Z_{k+1}+C_k^T\mathbb{E}_k(Z_{k+1}w_k),\\[1mm]
Z_N=GX_N,~~k\in \mathbb{T}_t
\end{array}
\right.
\end{eqnarray*}
with $X$ is given in (\ref{system-1}).

\end{corollary}

\emph{Proof.} From Proposition \ref{proposition-operator-Hilbert}, we need only to derive the expression of $\Theta_1^tu+\Theta_2^tx$. In fact, from (\ref{operator-theta-1})-(\ref{adjoint-V-2}) we have
\begin{eqnarray*}
(\Theta_1^tu+\Theta_2^tx)_k=
R_ku_k+B_k^T\mathbb{E}_{k-d}(V^1_{k+1}+V^2_{k+1})+D_k^T\mathbb{E}_{k-d}\big{(}(V^1_{k+1}+V^2_{k+1})w_k\big{)},~k\in \mathbb{T}_{t}.
\end{eqnarray*}
Noting (\ref{system-operator}) and (\ref{system-operator-2}), we will have (\ref{Lemma-difference-0}).\endpf 

\section{Problem (LQ) for a fixed time-state initial pair}\label{section--fixed}

In this section, we will study Problem (LQ) for the fixed initial pair $(t,x)$, and the general case of Problem (LQ) for all the initial pairs will be introduced and studied in the next section.
Throughout this section, Problem (LQ) for the fixed initial pair $(t,x)$ will be simply denoted as Problem (LQ)$_{tx}$.
%

Throughout this paper, $\mathbb{E}_{k-d}$ is understood as $\mathbb{E}_t$ if $k\in \{t,...,t+d-1\}$ (i.e., $k-d<t$). From Proposition \ref{proposition-operator-Hilbert} and Corollary \ref{Lemma-difference}, the following theorem is straightforward.

\begin{theorem}\label{Theorem--Nece-Suff-fixed}
The following statements are equivalent.

(i) Problem (LQ)$_{tx}$ is solvable.

(ii) The following assertions hold.
\begin{itemize}
\item[]
\begin{itemize}
\item[a)] There exists a $u^{{t,x,*}}\in \mathcal{U}_{ad}^t$ such that the stationary condition
\begin{eqnarray}\label{stationary-condition}
R_k u^{t,x,*}_k+B_k^T\mathbb{E}_{k-d}Z^{t,x,*}_{k+1}+D_k^T\mathbb{E}_{k-d}(Z^{t,x,*}_{k+1}w_k)=0, ~~a.s., ~~k\in \mathbb{T}_{t}
\end{eqnarray}
is satisfied, where $Z^{t,x,*}$ is the backward state of the following FBS$\Delta$E
\begin{eqnarray}\label{FBSDE}
\left\{\begin{array}{l}
X^{t,x,*}_{k+1}=\big{(}A_{k}X^{t,x,*}_k+B_{k}u^{t,x,*}_{k}\big{)}+\big{(}C_{k}X^{t,x,*}_k+D_{k}u^{t,x,*}_{k}\big{)}w_k, \\[1mm]
Z^{t,x,*}_k=Q_k X^{t,x,*}_k+A_k^T\mathbb{E}_k Z^{t,x,*}_{k+1}+C_k^T\mathbb{E}_k(Z^{t,x,*}_{k+1}w_k),\\[1mm]
X^{t,x,*}_t=x,~~Z^{t,x,*}_N=GX^{t,x,*}_N,~~k\in \mathbb{T}_t.
\end{array}
\right.
\end{eqnarray}

\item[b)] The convexity condition
\begin{eqnarray}\label{convexity-condition}
\inf_{{u}\in \mathcal{U}_{ad}^t}J(t,0;{u})\geq 0
\end{eqnarray}
holds.
\end{itemize}
\end{itemize}

Under any of the above conditions, $u^{{t,x,*}}$ in (ii) is an optimal control of Problem (LQ)$_{tx}$.

\end{theorem}

%
%
%
%


\subsection{Stationary condition}

%
%
%
By the stationary condition (\ref{stationary-condition}) and a backward procedure of calculations, we can get the following discrete-time Riccati-like equations %
%
\begin{eqnarray}\label{Riccati-1}
\left\{
\begin{array}{l}
P^{(0)}_k=Q_k+A_k^T\big{(}P^{(0)}_{k+1}+P^{(1)}_{k+1}\big{)}A_k+C_k^TP^{(0)}_{k+1}C_k,\\[1mm]
P^{(i)}_k=A_k^TP^{(i+1)}_{k+1} A_k,~~i=1,...,d-1,\\[1mm]
P^{(d)}_k=-H_k^TW_k^\dagger H_k,\\[1mm]
P^{(0)}_N=G,~~P^{(j)}_N=0,~~j=1,...,d,\\[1mm]
k\in \mathbb{T}_{t+d}=\{t+d,...,N-1\},
\end{array}
\right.
\end{eqnarray}
\begin{eqnarray}\label{Riccati-2}
\left\{
\begin{array}{l}
P^{(0)}_k=Q_k+A_k^T\big{(}P^{(0)}_{k+1}+P^{(1)}_{k+1}\big{)}A_k+C_k^TP_{k+1}^{(0)}C_k,\\[1mm]
P^{(i)}_{k}=A_k^TP^{(i+1)}_{k+1}A_k,~~i=1,...,k-t-1,\\[1mm]
P^{(k-t)}_k=A_k^TP^{(k+1-t)}_{k+1}A_k-H_k^TW_k^\dagger H_k,\\[1mm]
k\in \{t+2,...,t+d-1\},
\end{array}
\right.
\end{eqnarray}
and
\begin{eqnarray}\label{Riccati-3}
\left\{
\begin{array}{l}
\left\{
\begin{array}{l}
P^{(0)}_{t+1}=Q_{t+1}+A_{t+1}^T\big{(}P^{(0)}_{t+2}+P^{(1)}_{t+2}\big{)}A_{t+1}+C_{t+1}^TP_{t+2}^{(0)}C_{t+1},\\[1mm]
%
%
P^{(1)}_{t+1}=A_{t+1}^TP^{(2)}_{{t+2}}A_{t+1}-H_{t+1}^TW_{t+1}^\dagger H_{t+1},
\end{array}
\right.\\[2mm]
P^{(0)}_t=Q_{t}+A_{t}^T\big{(}P^{(0)}_{t+1}+P^{(1)}_{t+1}\big{)}A_{t}+C_{t}^TP_{t+1}^{(0)}C_{t}-H_{t}^TW_{t}^\dagger H_{t},
\end{array}
\right.
\end{eqnarray}
where
\begin{eqnarray}\label{W}
W_{k}=\ds\left\{
\begin{array}{ll}
R_k+\sum_{i=0}^{d}B_k^T P^{(i)}_{k+1} B_k+D_k^TP^{(0)}_{k+1} D_k,&~~k\in \mathbb{T}_{t+d},\\[2mm]
R_k+\sum_{i=0}^{k+1-t}B_k^T P^{(i)}_{k+1} B_k+D_k^TP^{(0)}_{k+1} D_k,&~~k\in \{t,...,t+d-1\},
\end{array}
\right.
\end{eqnarray}
and
\begin{eqnarray}\label{H}
H_{k}=\ds\left\{
\begin{array}{ll}
\sum_{i=0}^{d}B_k^T P^{(i)}_{k+1}A_k+D_k^TP^{(0)}_{k+1}C_k,&~~k\in \mathbb{T}_{t+d},\\[2mm]
\sum_{i=0}^{k+1-t}B_k^T P^{(i)}_{k+1}A_k+D_k^TP^{(0)}_{k+1}C_k,&~~k\in \{t,...,t+d-1\}.
\end{array}
\right.
\end{eqnarray}
%
%
%
%
%
Furthermore, the backward state of FBS$\Delta$E (\ref{FBSDE}) can be expressed via the forward state and the solution of (\ref{Riccati-1})-(\ref{Riccati-3}). Due to the $d$-step-lagged information structure, the Riccati-like equations are much more complicated than the standard discrete-time Riccati equation.
In fact, we have the following equivalent characterization of the stationary condition.

\begin{theorem}\label{Theorem-necessary}
The following statements are equivalent.

(i) The stationary condition of (\ref{stationary-condition}) is satisfied for some $u^{t,x,*}\in \mathcal{U}_{ad}^t$.

(ii) The following condition
\begin{eqnarray}\label{Theorem-necessary-equality}
H_k\mathbb{E}_{k-d}X_k^{t,x,*}\in \mbox{Ran}(W_k),~~a.s.,~~k\in \mathbb{T}_{t}
\end{eqnarray}
is satisfied, where $W_k, H_k, k\in \mathbb{T}_t$, are given in (\ref{W}) and (\ref{H}), and $X^{t,x,*}$ is given by the forward S$\Delta$E of
\begin{eqnarray}\label{FBSDE-2}
\left\{\begin{array}{l}
X^{t,x,*}_{k+1}=\big{(}A_{k}X^{t,x,*}_k+B_{k}u^{t,x,*}_{k}\big{)}+\big{(}C_{k}X^{t,x,*}_k+D_{k}u^{t,x,*}_{k}\big{)}w_k, \\[1mm]
Z^{t,x,*}_k=Q_k X^{t,x,*}_k+A_k^T\mathbb{E}_k Z^{t,x,*}_{k+1}+C_k^T\mathbb{E}_k(Z^{t,x,*}_{k+1}w_k),\\[1mm]
X^{t,x,*}_t=x,~~Z^{t,x,*}_N=GX^{t,x,*}_N,~~k\in \mathbb{T}_t
\end{array}
\right.
\end{eqnarray}
with
\begin{eqnarray}\label{Theorem-necessary-u}
u_k^{t,x,*}=
-W_k^\dagger H_k\mathbb{E}_{k-d}X^{t,x,*}_k,&~~k\in \mathbb{T}_{t}.
\end{eqnarray}
Furthermore, the backward state $Z^{t,x,*}$ of (\ref{FBSDE-2}) has the following expression
\begin{eqnarray}\label{Z-X}
Z_k^{t,x,*}=\left\{
\begin{array}{ll}
P_k^{(0)}X^{t,x,*}_k+P_k^{(1)}\mathbb{E}_{k-1}X^{t,x,*}_k+\cdots+P^{(k-t)}_{k}\mathbb{E}_tX_k^{t,x,*},&~~k\in \{t,...,t+d-1\}, \\[2mm]
P_k^{(0)}X^{t,x,*}_k+P_k^{(1)}\mathbb{E}_{k-1}X^{t,x,*}_k+\cdots+P^{(d)}_{k}\mathbb{E}_{k-d}X_k^{t,x,*},&~~k\in \mathbb{T}_{t+d},
\end{array}
\right.
\end{eqnarray}
where $P^{(i)}, i=0,...,d$, are given in (\ref{Riccati-1})-(\ref{Riccati-3}).


\end{theorem}

\emph{Proof.} See Appendix A. \endpf 

\begin{remark}
If $x$ in (\ref{FBSDE-2}) is 0, then $X^{t,0,*}_k=0, k\in \mathbb{T}_t$. In this case, the condition (\ref{Theorem-necessary-equality}) is naturally satisfied.

\end{remark}

\begin{remark}
From the proof of Theorem \ref{Theorem-necessary}, we know that the key technique is to decouple the FBS$\Delta$E (\ref{FBSDE}) by virtue of (\ref{stationary-condition}), i.e., find the expression (\ref{Z-X}). Due to the delayed information structure, at $k\in \{t,...,t+d\}$ the decision information set remains $\mathcal{F}_t$. For $k\in \{t,...,t+d-1\}$, $Z^{t,x,*}_k$ is a linear function of $X^{t,x,*}_k, \mathbb{E}_{k-1}X^{t,x,*}_k, \cdots, P^{(k-t)}_{k}\mathbb{E}_tX_k^{t,x,*}$, which differs from the case of $k\in \mathbb{T}_{t+d}$. This is why the Riccati-like equations are divided into several pieces (\ref{Riccati-1})-(\ref{Riccati-3}).
Letting $k=t, t+1$, then $k-t-1$ in (\ref{Riccati-2}) will be 0 and 1. Hence, (\ref{Riccati-3}) is not a special form of (\ref{Riccati-2}).

\end{remark}

\begin{remark}
Substituting (\ref{Theorem-necessary-u}) into the forward S$\Delta$E of (\ref{FBSDE-2}), we have
\begin{eqnarray}\label{closed-loop-equation}
\left\{
\begin{array}{l}
X^{t,x,*}_{k+1}=\big{(}A_{k}X^{t,x,*}_k-B_{k}W_k^\dagger H_k\mathbb{E}_{k-d}X^{t,x,*}_k\big{)}+\big{(}C_{k}X^{t,x,*}_k-D_{k}W_k^\dagger H_k\mathbb{E}_{k-d}X^{t,x,*}_k\big{)}w_k,\\[1mm]
X^{t,x,*}_t=x,~~k\in \mathbb{T}_{t}.
\end{array}
\right.
\end{eqnarray}
Taking conditional expectations in both sides of all the equations of (\ref{closed-loop-equation}), we have
\begin{eqnarray*}
\left\{
\begin{array}{l}
\mathbb{E}_{k+1-d}X^{t,x,*}_{k+1}=A_{k}\mathbb{E}_{k+1-d}X^{t,x,*}_k-B_{k}W_k^\dagger H_k\mathbb{E}_{k-d}X^{t,x,*}_k,\\[1mm]
\mathbb{E}_{t-d}X^{t,x,*}_t=x,~~k\in \mathbb{T}_{t}.
\end{array}
\right.
\end{eqnarray*}
For $k\in \mathbb{T}_{t+d}$ and by successive iterations, it holds
\begin{eqnarray*}\label{closed-loop-equation-2}
&&\mathbb{E}_{k+1-d}X^{t,x,*}_{k+1}=A_{k}A_{k-1}\cdots A_{k+1-d}X^{t,x,*}_{k+1-d}- B_{k}W_k^\dagger H_k\mathbb{E}_{k-d}X^{t,x,*}_k\nonumber \\[1mm]
&&\hphantom{\mathbb{E}_{k+1-d}X^{t,x,*}_{k+1}=}\ds-\sum_{i=0}^{d-2}A_{k}\cdots A_{k-i}B_{k-i-1}W_{k-i-1}^\dagger H_{k-i-1}\mathbb{E}_{(k-i-1)-d}X^{t,x,*}_{k-i-1},
\end{eqnarray*}
which is eventually a linear function of $X^{t,x,*}_{k+1-d},...,X^{t,x,*}_t$.
Similar expressions can be derived for the case of $k\in \{t,...,t+d-1\}$. Combining this and (\ref{closed-loop-equation}), we can get all the values of $\mathbb{E}_{k-d}X_k^{t,x,*}, k\in \mathbb{T}_t$. Hence, the control (\ref{Theorem-necessary-u}) can be easily implemented.

\end{remark}

The following result shows that the solution of (\ref{Riccati-1})-(\ref{Riccati-3}) can be calculated through \emph{a} set of  Riccati-like equations.

\begin{proposition}
Let
\begin{eqnarray}\label{Riccati-N}
\left\{
\begin{array}{l}
\bar{P}^{(0)}_k=Q_k+A_k^T\big{(}\bar{P}^{(0)}_{k+1}+\bar{P}^{(1)}_{k+1}\big{)}A_k+C_k^T\bar{P}^{(0)}_{k+1}C_k,\\[1mm]
\bar{P}^{(i)}_k=A_k^T\bar{P}^{(i+1)}_{k+1} A_k,~~i=1,...,d-1,\\[1mm]
\bar{P}^{(d)}_k=-\bar{H}_k^T\bar{W}_k^\dagger \bar{H}_k,\\[1mm]
\bar{P}^{(0)}_N=G,~~\bar{P}^{(j)}_N=0,~~j=1,...,d,\\[1mm]
k\in \mathbb{T}_{t},
\end{array}
\right.
\end{eqnarray}
where
\begin{eqnarray*}
\left\{
\begin{array}{l}
\bar{W}_{k}=
R_k+\sum_{i=0}^{d}B_k^T \bar{P}^{(i)}_{k+1} B_k+D_k^T\bar{P}^{(0)}_{k+1} D_k,\\[1mm]
\bar{H}_k=\sum_{i=0}^{d}B_k^T \bar{P}^{(i)}_{k+1}A_k+D_k^T\bar{P}^{(0)}_{k+1}C_k,\\[1mm]
k\in \mathbb{T}_t.
\end{array}
\right.
\end{eqnarray*}
Then
for (\ref{Riccati-1})-(\ref{Riccati-3}) it holds that
\begin{eqnarray}
P_k^{(i)}=\left\{
\begin{array}{ll}
\bar{P}_k^{(i)},&~~~k\in \mathbb{T}_{t+d},~i=0,...,d,\\[1mm]
\bar{P}_k^{(i)},&~~~k\in \{t+2,...,t+d-1\},~i=0,...,k-t-1,\\[1mm]
\bar{P}^{(k-t)}_k+\cdots+\bar{P}^{(d)}_{k},&~~~k\in \{t,...,t+d-1\},~i=k-t.\\[1mm]
\end{array}
\right.
\end{eqnarray}

\end{proposition}

\emph{Proof}. $P_{k}^{(i)}=\bar{P}_k^{(i)}$ follows from their expressions for the case with $k\in \mathbb{T}_{t+d}, i=0,...,d$ and the case with $k\in \{t+2,...,t+d-1\}, i=0,...,k-t-1$. For $k=t+d-1$ and $i=d-1$,
\begin{eqnarray*}
&&P^{(d-1)}_{t+d-1}=A_{t+d-1}^TP^{(d)}_{t+d}A_{t+d-1}-H^T_{t+d-1}W^\dagger_{t+d-1}H_{t+d-1}\\[1mm]
&&\hphantom{P^{(d-1)}_{t+d-1}}=A_{t+d-1}^T\bar{P}^{(d)}_{t+d}A_{t+d-1}-\bar{H}^T_{t+d-1}\bar{W}^\dagger_{t+d-1}\bar{H}_{t+d-1}\\[1mm]
&&\hphantom{P^{(d-1)}_{t+d-1}}=\bar{P}^{(d-1)}_{t+d-1}+\bar{P}^{(d)}_{t+d-1}.
\end{eqnarray*}
Furthermore, we have
\begin{eqnarray*}
&&P^{(d-2)}_{t+d-2}=A_{t+d-2}^TP^{(d-1)}_{t+d-1}A_{t+d-2}-H^T_{t+d-2}W^\dagger_{t+d-2}H_{t+d-2}\\[1mm]
&&\hphantom{P^{(d-2)}_{t+d-2}}=A_{t+d-2}^T\big{(}\bar{P}^{(d-1)}_{t+d-1}+\bar{P}^{(d)}_{t+d-1}\big{)}A_{t+d-2}^T-\bar{H}^T_{t+d-2}\bar{W}^\dagger_{t+d-2}\bar{H}_{t+d-2}\\[1mm]
&&\hphantom{P^{(d-2)}_{t+d-2}}=\bar{P}^{(d-2)}_{t+d-2}+\bar{P}^{(d-1)}_{t+d-2}+\bar{P}^{(d)}_{t+d-2},
\end{eqnarray*}
where we have used the properties
\begin{eqnarray*}
&&H_{t+d-2}=\sum_{i=0}^{d-1}B_{t+d-2}^T P^{(i)}_{{t+d-1}}A_{t+d-2}+D_{t+d-2}^TP^{(0)}_{{t+d-1}}C_{t+d-2}\\[1mm]
&&\hphantom{H_{t+d-2}}=\sum_{i=0}^{d}B_{t+d-2}^T \bar{P}^{(i)}_{{t+d-1}}A_{t+d-2}+D_{t+d-2}^T\bar{P}^{(0)}_{{t+d-1}}C_{t+d-2}\\[1mm]
&&\hphantom{H_{t+d-2}}=\bar{H}_{t+d-2}
\end{eqnarray*}
and $W_{t+d-2}=\bar{W}_{t+d-2}$. By induction, we can achieve the conclusion. \endpf

\begin{remark}
(\ref{Riccati-N}) with $\bar{W}_k>0, k\in \mathbb{T}_t$, is first introduced in \cite{Zhang-huanshui2015}, which characterizes the unique solvability of stochastic LQ problem with input delay. We here will investigate Problem (LQ) with information transmission delay, and intend seeking more general conditions to ensure the solvability of Problem (LQ) for the case with a fixed initial pair and the case with all the initial pairs.

\end{remark}

\subsection{Convexity}

We now study the convexity condition. In what follows, the functional $u\mapsto J(t,x;u)$ is called convex if (\ref{convexity-condition}) holds.
%

%
%
%

\begin{lemma}\label{Lemma-Cost}
For any $u\in \mathcal{U}_{ad}^t$, it holds that
\begin{eqnarray}\label{add-subtr-1}
J(t,0;u)=\sum_{k=t}^{N-1}\mathbb{E}\Big{\{}(\mathbb{E}_{k-d}X^0_k)^TH_k^TW_k^{\dagger}H_k\mathbb{E}_{k-d}X^0_k+2(H_k\mathbb{E}_{k-d}X^0_k)^Tu_k+u_k^TW_ku_k \Big{\}}
\end{eqnarray}
with $X^0$ given by
\begin{eqnarray}\label{system-0}
\left\{\begin{array}{l}
X^0_{k+1}=\big{(}A_{k}X^0_k+B_{k}u_{k}\big{)}+\big{(}C_{k}X^0_k+D_{k}u_{k}\big{)}w_k, \\[1mm]
X^0_t=0,~~k\in \mathbb{T}_t.
\end{array}
\right.
\end{eqnarray}

\end{lemma}

\emph{Proof}. See Appendix B.  \endpf

As $W_k, k\in \mathbb{T}_t,$  are symmetric, there exist orthogonal matrices $F_k, k\in \mathbb{T}_t$, such that
\begin{eqnarray*}
W_k=F_k^T\left[
\begin{array}{cc}
\Sigma_{k}&0\\0&0
\end{array}
\right]F_k,~~k\in \mathbb{T}_t.
\end{eqnarray*}
In the above, $\Sigma_{k}\in \mathbb{R}^{r_k\times r_k}$ is a diagonal matrix with $r_k$ being the rank of $W_k$, whose diagonal elements are the nonzero eigenvalues of $W_k$. Hence, we have
\begin{eqnarray*}
W^\dagger_k=F_k^T\left[
\begin{array}{cc}
\Sigma^{-1}_{k}&0\\0&0
\end{array}
\right]F_k,~~k\in \mathbb{T}_t.
\end{eqnarray*}
Moreover, $F_k$ can be decomposed as $F_k^T=[(F^{(1)}_k)^T~~(F^{(2)}_k)^T]$, where the lines of $F^{(2)}_k\in \mathbb{R}^{(m-r_k)\times m}$ form a basis of $\mbox{Ker}(W_k)$ (the kernel of $W_k$).
Let
\begin{eqnarray*}
F_ku_k=\left[
\begin{array}{l}
F_k^{(1)}u_k\\[1mm]F_k^{(2)}u_k
\end{array}
\right],~~
L_k \triangleq F_kH_k=\left[
\begin{array}{l}
F_k^{(1)}H_k\\[1mm]F_k^{(2)}H_k
\end{array}
\right]\triangleq\left[
\begin{array}{l}
L_k^{(1)}\\[1mm]L^{(2)}_k
\end{array}
\right].
\end{eqnarray*}
Hence, (\ref{add-subtr}) becomes to
\begin{eqnarray}\label{Theorem-convex-0}
&&J(t,0;{u})%
=\sum_{k=t}^{N-1}\mathbb{E}\Big{[}\big{(}F_k^{(1)}u_k+\Sigma_k^{-1} L^{(1)}_k  \mathbb{E}_{k-d}X^0_k\big{)}^T \Sigma_k\big{(}F_k^{(1)}u_k+\Sigma_k^{-1} L^{(1)}_k  \mathbb{E}_{k-d}X^0_k\big{)}\Big{]}\nonumber\\
%
%
&&\hphantom{J(t,0;{u})=}+2\sum_{k=t}^{N-1}\mathbb{E}\Big{[}(L_k^{(2)}\mathbb{E}_{k-d}X^0_k)^TF_k^{(2)}u_k\Big{]}.
\end{eqnarray}
Note that the space spanned by lines of $F_k^{(1)}$ is $\mbox{Ran}(W_k)$, the range of $W_k$.
Let $\mathcal{U}_{ad}^t(\mbox{Ker})$ be a subset of $\mathcal{U}_{ad}^t$ such that for any $u\in \mathcal{U}_{ad}^t(\mbox{Ker})$, $u_k\in \mbox{Ker}(W_k), k\in \mathbb{T}_t$. Similarly, $\mathcal{U}_{ad}^t(\mbox{Ran})$ is defined as a subset of $\mathcal{U}_{ad}^t$ such that for any $u\in \mathcal{U}_{ad}^t(\mbox{Ran})$, $u_k\in \mbox{Ran}(W_k), k\in \mathbb{T}_t$.

By the above preparations, we have the following theorem, which gives necessary and sufficient conditions on the convexity of $u\mapsto J(t,x;u)$. To the best of our knowledge, it seems to be the first result to equivalently characterize the convexity of LQ problem.

\begin{theorem}\label{Theorem-convex-nece-suff}
The following statements are equivalent.

(i)  $u\mapsto J(t,x;u)$ is convex.

(ii) The following assertions hold.
\begin{itemize}
\item[]\begin{itemize}
\item[a)] The solution of Riccati-like equation set (\ref{Riccati-1})-(\ref{Riccati-3}) has the property $W_k\geq 0, k\in \mathbb{T}_t$.

\item[b)] For any $u\in \mathcal{U}_{ad}^t$, the condition
\begin{eqnarray}\label{Theorem-convex-nece-suff-WH}
H_k\mathbb{E}_{k-d}X^{0,u}_k\in \mbox{Ran}(W_k),~~a.s.,~~k\in \mathbb{T}_{t}
\end{eqnarray}
is satisfied, where $X^{0,u}$ is given by
\begin{eqnarray}\label{Theorem-convex-nece-suff-X}
\left\{\begin{array}{l}
X^{0,u}_{k+1}=\big{(}A_{k}X^{0,u}_k+B_{k}v^u_{k}\big{)}+\big{(}C_{k}X^{0,u}_k+D_{k}v^u_{k}\big{)}w_k, \\[1mm]
X^{0,u}_t=0,~~k\in \mathbb{T}_t
\end{array}
\right.
\end{eqnarray}
with
\begin{eqnarray}\label{Theorem-convex-nece-suff-v-u}
v^u_k=
u_k-W_k^\dagger H_k\mathbb{E}_{k-d}X_k^{0,u},~~k\in \mathbb{T}_{t}.
\end{eqnarray}

\end{itemize}

\end{itemize}
\end{theorem}

\emph{Proof}. $(i)\Rightarrow (ii)$.
As $u\mapsto J(t,x;u)$ is convex, from (\ref{Theorem-convex-0}) we have %
\begin{eqnarray*}
&&J(t,0;{u})%
=\sum_{k=t}^{N-1}\mathbb{E}\Big{[}\big{(}F_k^{(1)}u^{(1)}_k+\Sigma_k^{-1} L^{(1)}_k  \mathbb{E}_{k-d}X^0_k\big{)}^T \Sigma_k\big{(}F_k^{(1)}u^{(1)}_k+\Sigma_k^{-1} L^{(1)}_k  \mathbb{E}_{k-d}X^0_k\big{)}\Big{]}\nonumber\\
%
%
&&\hphantom{J(t,0;{u})=}+2\sum_{k=t}^{N-1}\mathbb{E}\Big{[}(L_k^{(2)}\mathbb{E}_{k-d}X^0_k)^TF_k^{(2)}u^{(2)}_k\Big{]} 
\nonumber\\
&&\hpm{J(t,0;{u})}\geq 0,
\end{eqnarray*}
where $u_k^{(1)}$ and $u_k^{(2)}$ are the projections of $u_k$ onto $\mbox{Ran}(W_k)$ and $\mbox{Ker}(W_k)$, respectively.
Then, it holds that
\begin{eqnarray}\label{Theorem-convex-1}
&&\inf_{u\in \mathcal{U}_{ad}^t(\mbox{Ran})}J(t,0;u)\nonumber \\
&&=\inf_{u\in\mathcal{U}_{ad}^t(\mbox{Ran})}\sum_{k=t}^{N-1}\mathbb{E}\Big{[}\big{(}F_k^{(1)}u_k+\Sigma_k^{-1} L^{(1)}_k  \mathbb{E}_{k-d}X^0_k\big{)}^T \Sigma_k\big{(}F_k^{(1)}u_k+\Sigma_k^{-1} L^{(1)}_k  \mathbb{E}_{k-d}X^0_k\big{)}\Big{]}\nonumber\\[1mm]
%
%
&&\geq 0.
\end{eqnarray}
Introduce a set
\begin{eqnarray*}
\widetilde{\mathcal{U}}^t_{ad}(\mbox{Ran})=\Big{\{}(F_t^{(1)}u_t,...,F_{N-1}^{(1)}u_{N-1})\,\big{|}\,u=\{u_t,...,u_{N-1}\}\in \mathcal{U}^t_{ad}(\mbox{Ran}) \Big{\}}.
\end{eqnarray*}
Note that for $k\in \mathbb{T}_t$, the lines, $\alpha_k^1,...,\alpha_k^{r_k}$, of $F^{(1)}_k$ form a basis of $W_k$.
For $u=\{u_t,...,u_{N-1}\}\in\mathcal{U}^t_{ad}(\mbox{Ran})$, for each $k\in \mathbb{T}_t$ there exist $\lambda_k^{1},...,\lambda_{k}^{r_k}\in \mathbb{R}$ such that
$
u_k=\sum_{i=1}^{r_k}\lambda_k^{i} (\alpha_k^i)^T
$.
Then,
\begin{eqnarray*}
F^{(1)}_ku_k=\sum_{i=1}^{r_k}\lambda_k^{i} \left[
\begin{array}{c}
\alpha_k^1\\ \vdots\\
\alpha_k^{r_k}
\end{array}
\right](\alpha_k^i)^T=\left[
\begin{array}{c}
\lambda_k^1\\ \vdots\\
\lambda_k^{r_k}
\end{array}
\right]\triangleq \lambda_k.
\end{eqnarray*}
For $k\in {T}_{t+d+1}=\{t+d+1,...,N-1\}$, $u_k$ is $\mathcal{F}_{k-d}$-measurable and $\mathbb{E}|u_k|^2< \infty$. This implies that $\lambda_k$ is $\mathcal{F}_{k-d}$-measurable and $\mathbb{E}|\lambda_k|^2< \infty$. Similar result holds for the case of $k\in \{t,...,t+d\}$.
Therefore,
\begin{eqnarray*}
\Big{\{}F_k^{(1)}u_k\,\big{|}\,u=\{u_t,...,u_{N-1}\}\in \mathcal{U}^t_{ad}(\mbox{Ran}) \Big{\}}=\left\{
\begin{array}{ll}
L^2_{\mathcal{F}}(t;\mathbb{R}^{r_k}),&~k\in \{t,...,t+d\},\\[1mm]
L^2_{\mathcal{F}}(k-d;\mathbb{R}^{r_k}),&~k\in \mathbb{T}_{t+d+1}.
\end{array}
\right.
\end{eqnarray*}
Here, $L^2_{\mathcal{F}}(t;\mathbb{R}^{r_k})$ and $L^2_{\mathcal{F}}(k-d;\mathbb{R}^{r_k}), k\in \mathbb{T}_{t+d}$, are similarly defined as (\ref{L^2-t}).
Therefore,  %
\begin{eqnarray}\label{U}
\hspace{-1mm}\widetilde{\mathcal{U}}^t_{ad}(\mbox{Ran})=L^2_{\mathcal{F}}(t;\mathbb{R}^{r_t})\times\cdots\times L^2_{\mathcal{F}}(t;\mathbb{R}^{r_{t+d}})\times L^2_{\mathcal{F}}(t+1;\mathbb{R}^{r_{t+d+1}})\times\cdots \times L^2_{\mathcal{F}}(N-1-d;\mathbb{R}^{r_{N-1}}).
\end{eqnarray}
Introduce a bounded linear operator from ${\mathcal{U}}^t_{ad}(\mbox{Ran})$ to $\widetilde{\mathcal{U}}^t_{ad}(\mbox{Ran})$
\begin{eqnarray*}\label{Theorem-convex-L}
\mathcal{L}u:~~(\mathcal{L}u)_k= F_k^{(1)}u_k+\Sigma_k^{-1} F^{(1)}_kH_k  \mathbb{E}_{k-d}X^0_k,~~ k\in \mathbb{T}_t.
\end{eqnarray*}
Here, $X^0$ is the solution of (\ref{system-0}).
We now prove that
$\mathcal{L}$ is a surjection.
In fact, for any $\theta\in \widetilde{\mathcal{U}}^t_{ad}(\mbox{Ran})$, let
\begin{eqnarray*}\label{Theorem-convex-6}
\left\{
\begin{array}{l}
\bar{X}^0_{k+1}=\big{\{}A_{k}\bar{X}^0_k+B_k(F_k^{(1)})^T\big{[}\theta_k-\Sigma_k^{-1} F^{(1)}_kH_k  \mathbb{E}_{k-d}\bar{X}^0_k\big{]}\big{\}}\\[1mm]
\hphantom{\bar{X}^0_{k+1}=}+\big{\{}C_{k}\bar{X}^0_k+D_k(F_k^{(1)})^T\big{[}\theta_k-\Sigma_k^{-1} F^{(1)}_kH_k  \mathbb{E}_{k-d}\bar{X}^0_k\big{]}\big{\}}w_k, \\[1mm]
\bar{X}^0_t=0,~~k\in \mathbb{T}_t
\end{array}
\right.
\end{eqnarray*}
and
\begin{eqnarray}\label{Theorem-convex-5}
u_k=(F_k^{(1)})^T\big{[}\theta_k-\Sigma_k^{-1} F^{(1)}_kH_k  \mathbb{E}_{k-d}\bar{X}^0_k\big{]},~~k\in \mathbb{T}_t.
\end{eqnarray}
%
%
%
%
%
%
%
Note that $u$ given in (\ref{Theorem-convex-5}) is in $\mathcal{U}_{ad}^t(\mbox{Ran})$.
As $F_k^{(1)}(F_k^{(1)})^T=I_{r_k}$, from  (\ref{Theorem-convex-5}) we have
\begin{eqnarray*}
\theta_k=(\mathcal{L}u)_k,~~k\in \mathbb{T}_t.
\end{eqnarray*}
Hence, $\mathcal{L}$ is a surjection defined from  ${\mathcal{U}}^t_{ad}(\mbox{Ran})$ to $\widetilde{\mathcal{U}}^t_{ad}(\mbox{Ran})$. From this, (\ref{Theorem-convex-1}) and the proof by contradiction, we have $\Sigma_k> 0, k\in \mathbb{T}_t$. This further implies $W_k\geq 0, k\in \mathbb{T}_t$. Then, a) is proved.


We now prove b).
Note that (\ref{Theorem-convex-nece-suff-v-u}) equals to
\begin{eqnarray*}
v^u_k=
u_k-(F_k^{(1)})^T\Sigma_k^{-1}F_k^{(1)}H_k\mathbb{E}_{k-d}X_k^{0,u},~~k\in \mathbb{T}_{t}.
\end{eqnarray*}
 For any $u\in \mathcal{U}_{ad}^t$, there exist $u^{(1)}\in \mathcal{U}_{ad}^t(\mbox{Ran})$ and $u^{(2)}\in \mathcal{U}_{ad}^t(\mbox{Ker})$ such that $u=u^{(1)}+u^{(2)}$, i.e., $u_k=u_k^{(1)}+u_k^{(2)}, k\in\mathbb{T}_t$. 
%
%
%
%
%
%
%
%
%
From (\ref{Theorem-convex-0}), for any $u\in \mathcal{U}^t_{ad}$ we have
\begin{eqnarray}\label{Theorem-convex-nece-suff-J}
&&J(t,0;v^{u})%
=\sum_{k=t}^{N-1}\mathbb{E}\Big{[}\big{(}F_k^{(1)}v^u_k+\Sigma_k^{-1} F^{(1)}_k H_k \mathbb{E}_{k-d}X^{0,u}_k\big{)}^T \Sigma_k\big{(}F_k^{(1)}v^u_k+\Sigma_k^{-1} F^{(1)}_kH_k  \mathbb{E}_{k-d}X^{0,u}_k\big{)}\Big{]}\nonumber\\
%
%
&&\hphantom{J(t,0;v^{u})=}+2\sum_{k=t}^{N-1}\mathbb{E}\Big{[}(F_k^{(2)}H_k\mathbb{E}_{k-d}X^{0,u}_k)^TF_k^{(2)}v^u_k\Big{]}\nonumber 
\\
&&\hphantom{J(t,0;v^{u})}=\sum_{k=t}^{N-1}\mathbb{E}\Big{[}\big{(}F_k^{(1)}u^{(1)}_k\big{)}^T \Sigma_k\big{(}F_k^{(1)}u^{(1)}_k\big{)}\Big{]}+2\sum_{k=t}^{N-1}\mathbb{E}\Big{[}\big{(}(F_k^{(2)})^TF_k^{(2)}H_k\mathbb{E}_{k-d}X^{0,u}_k\big{)}^Tu^{(2)}_k\Big{]}\nonumber\\
%
&&\hphantom{J(t,0;v^{u})}\geq 0.
\end{eqnarray}
In the above, we have used the properties $F^{(i)}_ku_k=F^{(i)}_ku^{(i)}_k, i=1,2$, and $F^{(2)}_k(F^{(1)}_k)^T=0$. From (\ref{Theorem-convex-nece-suff-J}), we must have
%
%
%
%
\begin{eqnarray}\label{Theorem-convex-nece-suff-WH-2}
(F_k^{(2)})^TF_k^{(2)}H_k\mathbb{E}_{k-d}X^{0,u}_k=0,~~a.s.,~~k\in \mathbb{T}_t.
\end{eqnarray}
Otherwise, we can select some $u$ such that $J(t,0;v^u)<0$. In fact, assume there exist $k_1\in \mathbb{T}_t$ and $\widehat{u}\in \mathcal{U}_{ad}^t$ such that
\begin{eqnarray*}
c_0=\mathbb{E}\big{|}(F_{k_1}^{(2)})^TF_{k_1}^{(2)}H_{k_1}\mathbb{E}_{{k_1}-d}X^{0,\widehat{u}}_{k_1}\big{|}^2> 0.
\end{eqnarray*}
Denote %
\begin{eqnarray*}
&&c_1=\sum_{k=t}^{N-1}\mathbb{E}\Big{[}\big{(}F_k^{(1)}\widehat{u}^{(1)}_k\big{)}^T \Sigma_k\big{(}F_k^{(1)}\widehat{u}^{(1)}_k\big{)}\Big{]},\\
%
%
&&c_2=2\sum_{k=t}^{k_1-1}\mathbb{E}\Big{[}\big{(}(F_k^{(2)})^TF_k^{(2)}H_k\mathbb{E}_{k-d}X^{0,\widehat{u}}_k\big{)}^T\widehat{u}^{(2)}_k\Big{]}.
\end{eqnarray*}
Introduce a new control
\begin{eqnarray*}
\widetilde{u}_{k}=
\left\{
\begin{array}{cl}
\widehat{u}^{(2)}_k,&k\in \{t,...,k_1-1\},\\[1mm]
-\frac{1+c_1+c_2}{2c_0}(F_{k_1}^{(2)})^TF_{k_1}^{(2)}H_{k_1}\mathbb{E}_{{k_1}-d}X^{0,\widehat{u}}_{k_1},& k=k_1,\\[1mm]
0,&k\in \{k_1+1,...,N-1\},
\end{array}
\right.
\end{eqnarray*}
which is in $\mathcal{U}^t_{ad}(\mbox{Ker})$.
Then, we have
\begin{eqnarray*}
&&J(t,0;v^{\widehat{u}^{(1)}+\widetilde{u}})=\sum_{k=t}^{N-1}\mathbb{E}\Big{[}\big{(}F_k^{(1)}\widehat{u}^{(1)}_k\big{)}^T \Sigma_k\big{(}F_k^{(1)}\widehat{u}^{(1)}_k\big{)}\Big{]}+2\sum_{k=t}^{k_1-1} \mathbb{E}\Big{[}\big{(}(F_k^{(2)})^TF_k^{(2)}H_k\mathbb{E}_{k-d}X^{0,\widehat{u}}_k\big{)}^T\widehat{u}^{(2)}_k\Big{]}\nonumber\\[1mm]
&&\hphantom{J(t,0;v^{\widehat{u}^{(1)}+\widetilde{u}})=}+2\mathbb{E}\Big{[}\big{(}(F_{k_1}^{(2)})^TF_{k_1}^{(2)}H_{k_1}\mathbb{E}_{{k_1}-d}X^{0,\widehat{u}}_{k_1}\big{)}^T{\widetilde{u}}_{k_1}\Big{]}\nonumber \\[1mm]
&&\hphantom{J(t,0;v^{\widehat{u}^{(1)}+\widetilde{u}})}=-\frac{1+c_1+c_2}{c_0}\mathbb{E}|(F_{k_1}^{(2)})^TF_{k_1}^{(2)}H_{k_1}\mathbb{E}_{{k_1}-d}X^{0,\widehat{u}}_{k_1}|^2+c_1+c_2\nonumber\\[1mm]
&&\hphantom{J(t,0;v^{\widehat{u}^{(1)}+\widetilde{u}})}= -1.
\end{eqnarray*}
This contradicts the convexity of $u\mapsto J(t,x;u)$. Hence, we have (\ref{Theorem-convex-nece-suff-WH-2}). By multiplying $F_k^{(2)}$ and noting $F^{(2)}_k(F^{(2)}_k)^T=I_{m-r_k}$, it holds that
\begin{eqnarray}\label{Theorem-convex-nece-suff-WH-3}
F_k^{(2)}H_k\mathbb{E}_{k-d}X^{0,u}_k=0,~~a.s.,~~k\in \mathbb{T}_t,
\end{eqnarray}
which is equivalent to (\ref{Theorem-convex-nece-suff-WH}).

$(ii)\Rightarrow (i)$. From the condition (ii), (\ref{Theorem-convex-0}), (\ref{Theorem-convex-nece-suff-WH}) and (\ref{Theorem-convex-nece-suff-WH-3}), we have for any $u\in \mathcal{U}_{ad}^t$
\begin{eqnarray}\label{Theorem-convex-nece-suff-J-2}
&&J(t,0;v^{u})=\sum_{k=t}^{N-1}\mathbb{E}\Big{[}\big{(}F_k^{(1)}u_k\big{)}^T \Sigma_k\big{(}F_k^{(1)}u_k\big{)}\Big{]}+2\sum_{k=t}^{N-1}\mathbb{E}\Big{[}\big{(}(F_k^{(2)})^TF_k^{(2)}H_k\mathbb{E}_{k-d}X^{0,u}_k\big{)}^Tu_k\Big{]}\nonumber \\[1mm]
&&\hphantom{J(t,0;v^{u})}=\sum_{k=t}^{N-1}\mathbb{E}\Big{[}\big{(}F_k^{(1)}u_k\big{)}^T \Sigma_k\big{(}F_k^{(1)}u_k\big{)}\Big{]}\nonumber\\[1mm]
&&\hphantom{J(t,0;v^{u})}\geq 0.
\end{eqnarray}
We now show
\begin{eqnarray}\label{Theorem-convex-nece-suff-U}
\Big{\{}v^u\,\big{|}\, u\in \mathcal{U}_{ad}^t \Big{\}}=\mathcal{U}_{ad}^t,
\end{eqnarray}
where $v^u$ is given by (\ref{Theorem-convex-nece-suff-v-u}). For any ${\widetilde{v}}\in \mathcal{U}_{ad}^t$, let
\begin{eqnarray}\label{Theorem-convex-nece-suff-u-1}
u_k={\widetilde{v}}_k+W_k^\dagger H_k\mathbb{E}_{k-d}\widetilde{X}_k^{0},~~k\in \mathbb{T}_{t},
\end{eqnarray}
where
\begin{eqnarray*}
\left\{
\begin{array}{l}
{\widetilde{X}}^0_{k+1}=\big{(}A_{k}{\widetilde{X}}^0_k+B_k{\widetilde{v}}_k\big{)}+\big{(}C_{k}{\widetilde{X}}^0_k+D_k{\widetilde{v}}_k\big{)}w_k, \\[1mm]
{\widetilde{X}}^0_t=0,~~k\in \mathbb{T}_t.
\end{array}
\right.
\end{eqnarray*}
We then have from (\ref{Theorem-convex-nece-suff-v-u}) (\ref{Theorem-convex-nece-suff-u-1}) that $v^{u}={\widetilde{v}}$. Hence, (\ref{Theorem-convex-nece-suff-U}) holds, which together with (\ref{Theorem-convex-nece-suff-J-2}) implies
\begin{eqnarray}\label{Theorem-convex-nece-suff-J-equa}
\inf_{u\in \mathcal{U}_{ad}^t}J(t,0;{u})=\inf_{u\in \mathcal{U}_{ad}^t}J(t,0;v^{u})\geq 0.
\end{eqnarray}
This completes the proof. \endpf

In the proof of $(ii)\Rightarrow(i)$ of Theorem \ref{Theorem-convex-nece-suff}, we have used a simple technique of control shifting ($u\mapsto v^u$).
To make it more clear, we state the following lemma, whose proof is omitted here.
Firstly, introduce a set; let $L^2(\mathbb{T}_t; \mathbb{R}^{m\times n})$ be a set of $\mathbb{R}^{m\times n}$-valued deterministic processes such that for any $\nu=\{\nu_k, k\in \mathbb{T}_t\}\in L^2(\mathbb{T}_t; \mathbb{R}^{m\times n})$, $\sum_{k=t}^{N-1}|\nu_k|^2<\infty$.

\begin{proposition}\label{Lemma-cost}
For any $\Phi\in L^2(\mathbb{T}_t; \mathbb{R}^{m\times n})$, the following assertions hold.

(i) The property
 $$\big{\{}{u}-\Phi\mathbb{E}_{\cdot-d}X \,\big{|}\, u\in \mathcal{U}_{ad}^t \big{\}}=\mathcal{U}_{ad}^t$$
is satisfied, where ${u}-\Phi\mathbb{E}_{\cdot-d}X$ is the control $\{{u}_k-\Phi_k\mathbb{E}_{k-d}X_k, ~~k\in \mathbb{T}_t\}
$
with
\begin{eqnarray*}
\left\{\begin{array}{l}
X_{k+1}=\big{(}A_{k}X_k-B_{k}\Phi_k\mathbb{E}_{k-d}X_k+B_ku_k\big{)}\\[1mm]
\hphantom{X^{u}_{k+1}=}+\big{(}C_{k}X_k-D_{k}\Phi_k \mathbb{E}_{k-d}X_k+D_ku_k\big{)}w_k, \\[1mm]
X_t=x,~~k\in \mathbb{T}_t.
\end{array}
\right.
\end{eqnarray*}

(ii) it holds that
\begin{eqnarray*}\label{}
\inf_{u\in \mathcal{U}_{ad}^t}J(t,x;{u})=\inf_{u\in \mathcal{U}_{ad}^t}J(t,x;{u}-\Phi\mathbb{E}_{\cdot-d}X).
\end{eqnarray*}

\end{proposition}

Furthermore, by some simple calculations, we can show that the convexity of $u\mapsto J(t,x;u)$ has a semi-global property in the sense of the following result.

\begin{proposition}
If $u\mapsto J(t,x;u)$ is convex, then for any $(k,\xi)\in \mathbb{T}_t\times \mathbb{R}^n$ and $u\in \mathcal{U}_{ad}^t$, $u|_{\mathbb{T}_k}\mapsto J(k,\xi;u|_{\mathbb{T}_k})$ is also convex. Here, $\mathbb{T}_k=\{k,...,N-1\}$, and $u|_{\mathbb{T}_k}$ is the restriction of $u$ on $\mathbb{T}_k$.

\end{proposition}

\subsection{The solvability of Problem (LQ)$_{tx}$}

Noting (\ref{Theorem-necessary-u}) and (\ref{FBSDE-2}), we have
\begin{eqnarray}\label{Theorem--Nece-Suff-fixed-sys}
\left\{\begin{array}{l}
X^{t,x,*}_{k+1}=\big{(}A_{k}X^{t,x,*}_k-B_{k}W_k^\dagger H_k \mathbb{E}_{k-d}X^{t,x,*}_k\big{)}+\big{(}C_{k}X^{t,x,*}_k-D_{k}W_k^\dagger H_k \mathbb{E}_{k-d}X^{t,x,*}_k\big{)}w_k, \\[1mm]
X^{t,x,*}_t=x,~~k\in \mathbb{T}_t
\end{array}
\right.
\end{eqnarray}
with property (\ref{Theorem-necessary-equality}). Letting $X^{x,u}=X^{t,x,*}+X^{0,u}$ and from (\ref{Theorem-convex-nece-suff-X}) (\ref{Theorem--Nece-Suff-fixed-sys}), we have
\begin{eqnarray}\label{Theorem--Nece-Suff-fixed-sys-u}
\left\{\begin{array}{l}
X^{x,u}_{k+1}=\big{(}A_{k}X^{x,u}_k-B_{k}W_k^\dagger H_k \mathbb{E}_{k-d}X^{x,u}_k+B_ku_k\big{)}\\[1mm]
\hphantom{X^{x,u}_{k+1}=}+\big{(}C_{k}X^{x,u}_k-D_{k}W_k^\dagger H_k \mathbb{E}_{k-d}X^{x,u}_k+D_ku_k\big{)}w_k, \\[1mm]
X^{x,u}_t=x,~~k\in \mathbb{T}_t.
\end{array}
\right.
\end{eqnarray}
From Theorem \ref{Theorem--Nece-Suff-fixed}, Theorem \ref{Theorem-necessary} and Theorem \ref{Theorem-convex-nece-suff}, we then have the following necessary and sufficient conditions on the existence of optimal control of Problem (LQ)$_{tx}$.

\begin{theorem}\label{Theorem--Nece-Suff-fixed-final}
The following statements are equivalent.

(i) Problem (LQ)$_{tx}$ is solvable.

(ii) The following assertions hold
\begin{itemize}
\item[]\begin{itemize}
\item[a)]
The solution of Riccati-like equation set (\ref{Riccati-1})-(\ref{Riccati-3}) has the property $W_k\geq 0, k\in \mathbb{T}_t$.

\item[b)] For any $u\in \mathcal{U}_{ad}^t$, the condition
 \begin{eqnarray}\label{Theorem--Nece-Suff-fixed-final-HW}
H_k\mathbb{E}_{k-d}X^{x,u}_k\in \mbox{Ran}(W_k),~~a.s.,~~k\in \mathbb{T}_{t}
\end{eqnarray}
is satisfied, where $X^{x,u}$ is the solution of (\ref{Theorem--Nece-Suff-fixed-sys-u}).
\end{itemize}
\end{itemize}

Under any of the above conditions, the following control
\begin{eqnarray}
u_k^{t,x,*}=-W_k^\dagger H_k\mathbb{E}_{k-d}X^{t,x,*}_k,&~~k\in \mathbb{T}_{t}
\end{eqnarray}
is an optimal control of Problem (LQ)$_{tx}$, where $X^{t,x,*}$ is given by (\ref{Theorem--Nece-Suff-fixed-sys}).

\end{theorem}

\emph{Proof}. The equivalence between (i) and (ii) follows from the construction of $X^{x,u}$. From Proposition \ref{Lemma-cost}, we have
\begin{eqnarray}\label{cost-equa}
\inf_{u\in \mathcal{U}_{ad}^t}J(t,x;{u})=\inf_{u\in \mathcal{U}_{ad}^t}J(t,x;v^{u})
\end{eqnarray}
with $v^u=\{v^u_k=u_k-W_k^\dagger H_k \mathbb{E}_{k-d}X^{x,u}_k, k\in \mathbb{T}_t\}$. Under any of (i) (ii) and similar to (\ref{add-subtr}), we have %
\begin{eqnarray}\label{cost-equa-2}
&&\hspace{-1em}J(t,x;v^u)=x^TP^{(0)}_tx+\sum_{k=t}^{N-1}\mathbb{E}\Big{\{}(\mathbb{E}_{k-d}X^{x,u}_k)^TH_k^TW_k^{\dagger}H_k\mathbb{E}_{k-d}X^{x,u}_k\nonumber \\
&&\hspace{-1em}\hphantom{J(t,x;v^u)=}+2(H_k\mathbb{E}_{k-d}X^{x,u}_k)^Tv^u_k +(v^u_k)^TW_kv^u_k \Big{\}}\nonumber \\
&&\hspace{-1em}\hphantom{J(t,x;v^u)}=x^TP^{(0)}_tx+\sum_{k=t}^{N-1}\mathbb{E}\Big{\{}u_k^TW_ku_k+u_k^T\big{(}H_k\mathbb{E}_{k-d}X^{x,u}_{k} -W_kW_k^\dagger H_k \mathbb{E}_{k-d}X^{x,u}_{k} \big{)} \Big{\}}\nonumber  \\
&&\hspace{-1em}\hphantom{J(t,x;v^u)}=x^TP^{(0)}_tx+\sum_{k=t}^{N-1}\mathbb{E}\big{(}u_k^TW_ku_k\big{)}\nonumber  \\
&&\hspace{-1em}\hphantom{J(t,x;v^u)}\geq x^TP^{(0)}_tx,
\end{eqnarray}
where for $u=0$ the equality holds. In the above, we have used the property (\ref{Theorem--Nece-Suff-fixed-final-HW}),
which is equivalent to
\begin{eqnarray*}
(I-W_{k} W_{k}^\dagger) H_{k}\mathbb{E}_{k-d}X^{x,u}_{k}=0,~~a.s., ~~k\in \mathbb{T}_t.
\end{eqnarray*}
By (\ref{cost-equa}) and (\ref{cost-equa-2}), we then achieve the conclusion. \endpf

Introduce a set
\begin{eqnarray*}
\mathcal{I}_t=\Big{\{}x\,\big{|}\,\mbox{Problem (LQ)}_{tx} \mbox{ is solvable}\Big{\}}.
\end{eqnarray*}

\begin{theorem}
$\mathcal{I}_t$ is either empty or a linear subspace of $\mbox{Ker}[(I-W_{t} W_{t}^\dagger) H_{t}]$.

\end{theorem}

\emph{Proof}. Letting  $u=0$ in (\ref{Theorem--Nece-Suff-fixed-sys-u}), we have $X^{x,0}=X^{t,x,*}$, which is given in (\ref{Theorem--Nece-Suff-fixed-sys}).
For $x\in \mathcal{I}_t\neq \varnothing$, $x$ will be in $\mbox{Ker}[(I-W_{t} W_{t}^\dagger) H_{t}]$.
Then, for $x, x'\in \mathcal{I}_t, \alpha, \beta\in \mathbb{R}$, we have
\begin{eqnarray*}
\left\{\begin{array}{l}
\alpha X^{x,0}_{k+1}+\beta X^{x',0}_{k+1}=\Big{[}A_{k}\big{(}\alpha X^{x,0}_{k}+\beta X^{x',0}_{k}\big{)}-B_{k}W_k^\dagger H_k \big{(}\alpha\mathbb{E}_{k-d} X^{x,0}_{k+1}+\beta \mathbb{E}_{k-d} X^{x',0}_{k}\big{)}\Big{]}\\[1mm]
\hphantom{\alpha X^{x,0}_{k+1}+\beta X^{x',0}_{k+1}=}+\Big{[}C_{k}\big{(}\alpha X^{x,0}_{k}+\beta X^{x',0}_{k}\big{)}-D_{k}W_k^\dagger H_k \big{(}\alpha\mathbb{E}_{k-d} X^{x,0}_{k+1}+\beta \mathbb{E}_{k-d} X^{x',0}_{k}\big{)}\Big{]}w_k, \\[1mm]
\alpha X^{x,0}_{t}+\beta X^{x',0}_{t}=\alpha x+\beta x',~~k\in \mathbb{T}_t.
\end{array}
\right.
\end{eqnarray*}
Hence, $\alpha X^{x,0}+\beta X^{x',0}=X^{\alpha x+\beta x', 0}$, which further implies
\begin{eqnarray*}
H_k\mathbb{E}_{k-d}X^{{\alpha x+\beta x'},0}_k=\alpha H_k\mathbb{E}_{k-d} X^{x,0}_{k+1}+\beta H_k\mathbb{E}_{k-d} X^{x',0}_{k}\in \mbox{Ran}(W_k),~~k\in \mathbb{T}_{t}.
\end{eqnarray*}
Combining with the convexity, we know that Problem (LQ) is solvable at $(t, \alpha x+\beta x')$.  \endpf

To end this subsection, a sufficient condition is presented to ensure (\ref{Theorem--Nece-Suff-fixed-final-HW}).

\begin{theorem}\label{Theorem-final-fixed}
If $\mbox{Ran}(H_k)\subset \mbox{Ran}(W_k)$ (i.e., $W_{k} W_{k}^\dagger H_{k}=H_{k}$), $k\in \mathbb{T}_t$, then the condition   (\ref{Theorem--Nece-Suff-fixed-final-HW}) is satisfied.

\end{theorem}

\emph{Proof}. The proof is straightforward and hence, omitted here. \endpf

Combining the condition in Theorem \ref{Theorem-final-fixed} with a) of Theorem \ref{Theorem--Nece-Suff-fixed-final}, we can obtain in the next section much neater results of
Problem (LQ) (for all the initial pairs).

\subsection{The delay-free case }

Let us revisit the standard discrete-time stochastic LQ problem without time delay.

\textbf{Problem (LQ)$^{df}_{tx}$}. \emph{For the initial pair $(t,x)\in \mathbb{T}\times \mathbb{R}^n$, find a
$\bar{u}\in \mathcal{U}_{ad}^t$ such that
\begin{eqnarray}\label{Problem-LQ-df}
J(t,x;\bar{u}) = \inf_{u\in L^2_\mathcal{F}(\mathbb{T}_t; \mathbb{R}^m)}J(t,x;u).
\end{eqnarray}
In (\ref{Problem-LQ-df}), $J(t,x;u)$ is given in (\ref{cost-1}) and is subject to (\ref{system-1}), and
\begin{eqnarray*}
L^2_\mathcal{F}(\mathbb{T}_t; \mathbb{R}^m)=\Big{\{}\nu=\{\nu_k, k\in \mathbb{T}_{t}\}\,\big{|}\,\nu_k\mbox{ is }\mathcal{F}_k\mbox{-measurable}, k\in\mathbb{T}_{t}, \mbox{and }\mathbb{E}|\nu_k|^2<\infty, k\in \mathbb{T}_t \Big{\}}.
\end{eqnarray*}
}

Introduce the discrete-time Riccati equation
\begin{eqnarray}\label{Riccati}
\left\{
\begin{array}{l}
P_k=Q_{k}+A_{k}^T P_{k+1}A_{k}+C_{k}^TP_{k+1}C_{k}-H_{k}^TW_{k}^\dagger H_{k},\\[1mm]
P_N=G,~~k\in \mathbb{T}_t,
\end{array}
\right.
\end{eqnarray}
where
\begin{eqnarray*}
\left\{
\begin{array}{l}
W_{k}=R_k+B_k^T P_{k+1} B_k+D_k^TP_{k+1} D_k,\\[1mm]
H_{k}=B_k^T P_{k+1}A_k+D_k^TP_{k+1}C_k,\\[1mm]
k\in \mathbb{T}_t.
\end{array}
\right.
\end{eqnarray*}
Let $d=0$ in Theorem \ref{Theorem--Nece-Suff-fixed-final}, we have the following result.

\begin{theorem}\label{Theorem--Nece-Suff-fixed-final-stand}
The following statements are equivalent.

(i) Problem (LQ)$^{df}_{tx}$ is solvable.

(ii) The following assertions hold
\begin{itemize}
\item[]\begin{itemize}
\item[a)]
The solution of Riccati equation (\ref{Riccati}) has the property $W_k\geq 0, k\in \mathbb{T}_t$.

\item[b)] For any $u\in L^2_\mathcal{F}(\mathbb{T}_t; \mathbb{R}^m)$, the condition
 \begin{eqnarray}\label{Theorem-convex-nece-suff-WH-2-final}
H_kX^{x,u}_k\in \mbox{Ran}(W_k),~~a.s.,~~k\in \mathbb{T}_{t}
\end{eqnarray}
is satisfied, where $X^{x,u}$ is the solution of the following S$\Delta$E
\begin{eqnarray}\label{Theorem-convex-nece-suff-X-df}
\left\{\begin{array}{l}
X^{x,u}_{k+1}=\big{(}\bar{A}_{k}X^{x,u}_k+B_{k}u_{k}\big{)}+\big{(}\bar{C}_{k} X^{x,u}_k+D_{k}u_{k}\big{)}w_k, \\[1mm]
X^{x,u}_t=x,~~k\in \mathbb{T}_t
\end{array}
\right.
\end{eqnarray}
with
\begin{eqnarray*}
\bar{A}_k=A_{k}-B_kW_k^\dagger H_k,~~\bar{C}_k=C_{k}-D_kW_k^\dagger H_k,~~k\in \mathbb{T}_t.
\end{eqnarray*}

\end{itemize}
\end{itemize}
\end{theorem}

Firstly, let us take some observation. Let $V_k=(I-W_k^\dagger W_k)H_k, k\in \mathbb{T}_t$. Then, the condition (\ref{Theorem-convex-nece-suff-WH-2-final}) is equivalent to
\begin{eqnarray*}
V_kX^{x,u}_k=0,~~a.s.,~~k\in \mathbb{T}_t.
\end{eqnarray*}
Hence, at $k$, the attainable set of the system (\ref{Theorem-convex-nece-suff-X-df}) is a subset of $\mbox{Ker}(V_k)$.
Similarly, (\ref{Theorem--Nece-Suff-fixed-final-HW}) is relating to the property of the attainable set of system (\ref{Theorem--Nece-Suff-fixed-sys-u}).
To get conditions of (\ref{Theorem--Nece-Suff-fixed-final-HW}) (\ref{Theorem-convex-nece-suff-WH-2-final})  that are easier to validated, we should in the future to study the attainable set of (\ref{Theorem-convex-nece-suff-X-df}) and (\ref{Theorem--Nece-Suff-fixed-sys-u}), which is further related to the controllability of linear S$\Delta$Es.


Letting the initial pair $(t,x)$ vary in the product space $\mathbb{T}\times \mathbb{R}^n$, we get a family of LQ problems; from Theorem \ref{Theorem--Nece-Suff-fixed-final-stand}, we have an equivalent characterization of the solvability of this family of LQ problems.

\begin{proposition}\label{Proposition-delay-free}
The following statements are equivalent.

(i) For any $(t,x)\in \mathbb{T}\times \mathbb{R}^n$, Problem (LQ)$^{df}_{tx}$ is solvable.

(ii) The constrained Riccati equation
\begin{eqnarray}\label{Riccati-0}
\left\{
\begin{array}{l}
P_k=Q_{k}+A_{k}^T P_{k+1}A_{k}+C_{k}^TP_{k+1}C_{k}-H_{k}^TW_{k}^\dagger H_{k},\\[1mm]
P_N=G,\\[1mm]
W_k^\dagger W_kH_k=H_k,~W_k\geq 0,\\
k\in \mathbb{T}
\end{array}
\right.
\end{eqnarray}
is solvable in the sense $W_k^\dagger W_kH_k=H_k,~W_k\geq 0,~k\in \mathbb{T}_t$, where
\begin{eqnarray*}
\left\{
\begin{array}{l}
W_{k}=R_k+B_k^T P_{k+1} B_k+D_k^TP_{k+1} D_k,\\[1mm]
H_{k}=B_k^T P_{k+1}A_k+D_k^TP_{k+1}C_k,\\[1mm]
k\in \mathbb{T}.
\end{array}
\right.
\end{eqnarray*}

\end{proposition}

\emph{Proof}. For any $t\in \mathbb{T}$ and letting $k=t$ in (\ref{Theorem-convex-nece-suff-WH-2-final}), we have
$H_tx\in \mbox{Ran}(W_t)$, which holds for any $x\in \mathbb{R}^n$; equivalently, we have $\mbox{Ran}(H_t) \subset \mbox{Ran}(W_t)$ by considering the cases $x=e_i, i=1,...,n$. Here, $e_i$ is the $n$-dimensional vector, whose $i$-th entry is 1 and other entries are all 0. Combining this fact and Theorem \ref{Theorem--Nece-Suff-fixed-final-stand}, we then achieve the result.   \endpf

\begin{remark}
Proposition \ref{Proposition-delay-free} is a main result of \cite{Ait-Chen-Zhou-2002}, which presents a necessary and sufficient condition on the solvability of a family of LQ problems, that is,
\begin{center}
\{Problem (LQ)$^{df}_{tx}$ is solvable, for any $(t,x)\in \mathbb{T}\times \mathbb{R}^n$\}$\Leftrightarrow$ (\ref{Riccati-0}) is solvable.
\end{center}
In contrast, Theorem \ref{Theorem--Nece-Suff-fixed-final-stand} just characterizes the solvability of Problem (LQ)$^{df}_{tx}$. The proof of Proposition  \ref{Proposition-delay-free} shows that Theorem  \ref{Theorem--Nece-Suff-fixed-final-stand} implies Proposition  \ref{Proposition-delay-free}. However, the equivalence between (i) and (ii) of Theorem  \ref{Theorem--Nece-Suff-fixed-final-stand} cannot be proved by virtue of Proposition  \ref{Proposition-delay-free}.
Hence, Theorem \ref{Theorem--Nece-Suff-fixed-final-stand} is a new result even for standard LQ problems; a key step to access to Theorem \ref{Theorem--Nece-Suff-fixed-final-stand} is that we have obtained (for the first time) an equivalent characterization of the convexity of the cost functional.

\end{remark}


\section{Problem (LQ) for all the time-state initial pairs} \label{section--all}

\subsection{The solvability of Problem (LQ)}

In this section, we will study Problem (LQ) for all the initial pairs.
To begin with, we introduce versions of Problem (LQ) (for the initial pair $(t,x)$). If $k\in \{t,...,t+d-1\}$, let
\begin{eqnarray}\label{admissible-control--k-1}
\mathcal{U}_{ad}^k=\Big{\{}u=\{u_k,u_{k+1},...,u_{N-1}\}\,\big{|}\,u\in \big{(}l^2_\mathcal{F}(t; \mathbb{R}^n)\big{)}^{t+d-k}\times l^2_\mathcal{F}(\mathbb{T}_{t}^{-d}; \mathbb{R}^m)  \Big{\}};
\end{eqnarray}
if $k\in \mathbb{T}_{t+d}$, let
\begin{eqnarray}\label{admissible-control--k-2}
\mathcal{U}_{ad}^k=\Big{\{}u=\{u_k,u_{k+1},...,u_{N-1}\}\,\big{|}\,u\in l^2_\mathcal{F}(\mathbb{T}_{k-d}^{-d}; \mathbb{R}^m)  \Big{\}}.
\end{eqnarray}
In (\ref{admissible-control--k-2}), $l^2_\mathcal{F}(\mathbb{T}^{-d}_{k-d}; \mathbb{R}^m)$ is a set of
$\mathbb{R}^m$-valued processes with $\mathbb{T}_{k-d}^{-d}=\{k-d,...,N-1-d\}$ such that for any its element $\nu=\{\nu_\ell, \ell\in \mathbb{T}_{k-d}^{-d}\}$, $\nu_\ell$ is $\mathcal{F}_{\ell}$-measurable and
$\sum_{\ell=k-d}^{N-1-d}\mathbb{E}|\nu_\ell|^2<\infty$.

Starting from the initial pair $(k,\xi)\in \mathbb{T}_t\times \mathbb{R}^n$, (\ref{system-1}) (\ref{cost-1}) become to
\begin{eqnarray}\label{system-k}
\left\{\begin{array}{l}
X_{\ell+1}=\big{(}A_{\ell}X_\ell+B_{\ell}u_{\ell}\big{)}+\big{(}C_{\ell}X_\ell+D_{\ell}u_{\ell}\big{)}w_\ell, \\[1mm]
X_k=\xi,~~\ell\in \mathbb{T}_k=\{k,...,N-1\},
\end{array}
\right.
\end{eqnarray}
and
\begin{eqnarray}\label{cost-k}
&&\hspace{-2em}J(k,\xi; u)
 =\sum_{\ell=k}^{N-1}\mathbb{E}\big{[}X_\ell^TQ_{\ell}X_\ell+u_{\ell}^TR_{\ell}u_{\ell}\big{]}+\mathbb{E}\big{[}X_N^TGX_N\big{]}.
\end{eqnarray}
Problem (LQ) for the initial pair $(k,\xi)$ is referred to as the case that  minimizes (\ref{cost-k}) over $\mathcal{U}^k_{ad}$ (subject to (\ref{system-k})).

\begin{definition}
(i) Problem (LQ) is said to be finite at $(k,\xi)\in \mathbb{T}_t\times \mathbb{R}^n$, if
\begin{eqnarray}\label{defi--finite}
\inf_{u\in \mathcal{U}^k_{ad}}J(k,\xi; u)>-\infty.
\end{eqnarray}
Problem (LQ) is said to be finite if (\ref{defi--finite}) holds for any initial pair $(k,\xi)\in \mathbb{T}_t\times \mathbb{R}^n$.

(ii) Problem (LQ) is said to be (uniquely) solvable at $(k,\xi)\in \mathbb{T}_t\times \mathbb{R}^n$ if there exists a (unique) $\bar{u}\in \mathcal{U}^k_{ad}$
such that
\begin{eqnarray*}
J(k,\xi;\bar{u}) = \inf_{u\in \mathcal{U}_{ad}^k}J(k,\xi;u).
\end{eqnarray*}
In this case, $\bar{u}$ is called an optimal control of Problem (LQ) for the initial pair $(k,\xi)$. Problem (LQ) is said to be (uniquely) solvable if it is solvable at any initial pair $(k,\xi)\in \mathbb{T}_t \times \mathbb{R}^n$.
\end{definition}

To study the finiteness of Problem (LQ), introduce the following coupled LMEIs  (\ref{Riccati-inequality-1})-(\ref{Riccati-inequality-3})
\begin{eqnarray}\label{Riccati-inequality-1}
\left\{
\begin{array}{l}
P^{(0)}_k\leq Q_k+A_k^T\big{(}P^{(0)}_{k+1}+P^{(1)}_{k+1}\big{)}A_k+C_k^TP^{(0)}_{k+1}C_k,\\[1mm]
P^{(i)}_k=A_k^TP^{(i+1)}_{k+1} A_k,~~i=1,...,d-1,\\[1mm]
\left[
\begin{array}{cc}
-P^{(d)}_k&H_k^T\\H_k&W_k
\end{array} \right]\geq 0,\\[1mm]
P^{(0)}_N\leq G,~~P^{(j)}_N=0,~~j=1,...,d,\\[1mm]
k\in \mathbb{T}_{t+d},
\end{array}
\right.
\end{eqnarray}
\begin{eqnarray}\label{Riccati-inequality-2}
\left\{
\begin{array}{l}
P^{(0)}_k\leq Q_k+A_k^T\big{(}P^{(0)}_{k+1}+P^{(1)}_{k+1}\big{)}A_k+C_k^TP_{k+1}^{(0)}C_k,\\[1mm]
P^{(i)}_{k}=A_k^TP^{(i+1)}_{k+1}A_k,~~i=1,...,k-t-1,\\[1mm]
\left[
\begin{array}{cc}
A_k^TP^{(k+1-t)}_{k+1}A_k-P^{(k-t)}_k&H_k^T\\H_k&W_k
\end{array} \right]\geq 0,
\\[1mm]
k\in \{t+2,...,t+d-1\},
\end{array}
\right.
\end{eqnarray}
and
\begin{eqnarray}\label{Riccati-inequality-3}
\left\{
\begin{array}{l}
\left\{
\begin{array}{l}
P^{(0)}_{t+1}\leq Q_{t+1}+A_{t+1}^T\big{(}P^{(0)}_{t+2}+P^{(1)}_{t+2}\big{)}A_{t+1}+C_{t+1}^TP_{t+2}^{(0)}C_{t+1},\\[1mm]
%
%
\left[
\begin{array}{cc}
A_{t+1}^TP^{(2)}_{t+2}A_{t+1}-P^{(1)}_{t+1}&H_{t+1}^T\\H_{t+1}&W_{t+1}
\end{array} \right]\geq 0,%
\end{array}
\right.
\\[2mm]
\left[
\begin{array}{cc}
Q_{t}+A_{t}^T\big{(}P^{(0)}_{t+1}+P^{(1)}_{t+1}\big{)}A_{t}+C_{t}^TP_{t+1}^{(0)}C_{t}-P^{(0)}_{t}&H_{t}^T\\H_{t}&W_{t}
\end{array} \right]\geq 0,
%
%
%
\end{array}
\right.
\end{eqnarray}
where
\begin{eqnarray*}\label{Riccati-inequality-W}
W_{k}=\ds\left\{
\begin{array}{ll}
R_k+\sum_{i=0}^{d}B_k^T P^{(i)}_{k+1} B_k+D_k^TP^{(0)}_{k+1} D_k,&~~k\in \mathbb{T}_{t+d},\\[2mm]
R_k+\sum_{i=0}^{k+1-t}B_k^T P^{(i)}_{k+1} B_k+D_k^TP^{(0)}_{k+1} D_k,&~~k\in \{t,...,t+d-1\},
\end{array}
\right.
\end{eqnarray*}
and
\begin{eqnarray*}\label{Riccati-inequality-H}
H_{k}=\ds\left\{
\begin{array}{ll}
\sum_{i=0}^{d}B_k^T P^{(i)}_{k+1}A_k+D_k^TP^{(0)}_{k+1}C_k,&~~k\in \mathbb{T}_{t+d},\\[2mm]
\sum_{i=0}^{k+1-t}B_k^T P^{(i)}_{k+1}A_k+D_k^TP^{(0)}_{k+1}C_k,&~~k\in \{t,...,t+d-1\}.
\end{array}
\right.
\end{eqnarray*}

If exists, the solution of (\ref{Riccati-inequality-1})-(\ref{Riccati-inequality-3}) is denoted as $(P^{(0)},...,P^{(d)})$. Let
\begin{eqnarray}\label{M}
\mathcal{M}=\Big{\{}(P^{(0)},...,P^{(d)})\,\big{|}\,(P^{(0)},...,P^{(d)})\mbox{ is a solution of } (\ref{Riccati-inequality-1})\mbox{-}(\ref{Riccati-inequality-3})\Big{\}}.
\end{eqnarray}
Based on the solution of (\ref{Riccati-inequality-1})-(\ref{Riccati-inequality-3}), we have the following lemma, whose proof is similar to that of Lemma \ref{Lemma-Cost} and is omitted here.

\begin{lemma}\label{Lemma-add-subtr}
Let $(P^{(0)},...,P^{(d)})\in \mathcal{M}\neq \varnothing$. Then the following statements hold.

(i) For $k\in \mathbb{T}_{t+d}$, it holds that
\begin{eqnarray*}\label{add-subtr-2}
&&J(k,\xi;{u}) =\sum_{\ell=k}^{N-1}\mathbb{E}\Big{\{}X_\ell^T\big{[}Q_\ell+A_\ell^T\big{(}P^{(0)}_{\ell+1}+P^{(1)}_{\ell+1}\big{)}A_\ell^T+C_\ell^TP^{(0)}_{\ell+1}C_\ell -P_{\ell}^{(0)}\big{]}X_\ell\Big{\}}\nonumber \\
&&\hphantom{J(k,\xi;{u}) =}+\sum_{\ell=k}^{N-1}\mathbb{E}\left\{\left[\begin{array}{cc}
\mathbb{E}_{\ell-d}X_\ell\\
u_\ell
\end{array}\right]^T\left[
\begin{array}{cc}
-P^{(d)}_\ell &H_\ell^T\\H_\ell&W_\ell
\end{array} \right] \left[\begin{array}{cc}
\mathbb{E}_{\ell-d}X_\ell\\
u_\ell
\end{array}\right]\right\}\\
&&\hphantom{J(k,\xi;{u})=}
+\xi^T\Big{(}\sum_{i=0}^dP^{(i)}_k\Big{)}\xi+\mathbb{E}\big{[}X_N^T\big{(}G-P_N^{(0)}\big{)}X_N \big{]}\\
&&\hphantom{J(k,\xi;{u})}\geq \xi^T\Big{(}\sum_{i=0}^dP^{(i)}_k\Big{)}\xi.
\end{eqnarray*}

(ii) For $k\in \{t+1,...,t+d-1\}$, it holds that
\begin{eqnarray*}\label{add-subtr-7}
&&J(k,\xi;{u}) =\sum_{\ell=k}^{N-1}\mathbb{E}\Big{\{}X_\ell^T\big{[}Q_\ell+A_\ell^T\big{(}P^{(0)}_{\ell+1}+P^{(1)}_{\ell+1}\big{)}A_\ell^T+C_\ell^TP^{(0)}_{\ell+1}C_\ell -P_{\ell}^{(0)}\big{]}X_\ell\Big{\}}\nonumber \\
&&\hphantom{J(k,\xi;{u})=}+\sum_{\ell=t+d}^{N-1}\mathbb{E}\left\{\left[\begin{array}{cc}
\mathbb{E}_{\ell-d}X_\ell\\
u_\ell
\end{array}\right]^T\left[
\begin{array}{cc}
-P^{(d)}_\ell&H_\ell^T\\H_\ell&W_\ell
\end{array} \right] \left[\begin{array}{cc}
\mathbb{E}_{\ell-d}X_\ell\\
u_\ell
\end{array}\right]\right\}\nonumber \\
&&\hpm{J(k,\xi;{u})=}+\sum_{\ell=k}^{t+d-1}\mathbb{E}\left\{\left[\begin{array}{cc}
\mathbb{E}_{t}X_\ell\\
u_\ell
\end{array}\right]^T\left[
\begin{array}{cc}
A_\ell^TP^{(\ell+1-t)}_{\ell+1}A_\ell-P^{(\ell-t)}_\ell&H_\ell^T\\H_\ell&W_\ell
\end{array} \right]\left[\begin{array}{cc}
\mathbb{E}_{t}X_\ell\\
u_\ell
\end{array}\right]\right\}\nonumber  \\
&&\hphantom{J(k,\xi;{u})=}+\xi^T\Big{(}\sum_{i=0}^{k-t}P^{(i)}_k\Big{)}\xi+\mathbb{E}\big{[}X_N^T\big{(}G-P_N^{(0)}\big{)}X_N \big{]}\\
&&\hphantom{J(k,\xi;{u})}\geq\xi^T\Big{(}\sum_{i=0}^{k-t}P^{(i)}_k\Big{)}\xi.
%
\end{eqnarray*}

(iii) It holds that
\begin{eqnarray*}\label{add-subtr-3}
&&J(t,\xi;{u}) =\sum_{\ell=t+1}^{N-1}\mathbb{E}\Big{\{}X_\ell^T\big{[}Q_\ell+A_\ell^T\big{(}P^{(0)}_{\ell+1}+P^{(1)}_{\ell+1}\big{)}A_\ell^T+C_\ell^TP^{(0)}_{\ell+1}C_\ell -P_{\ell}^{(0)}\big{]}X_\ell\Big{\}}\nonumber \\
&&\hphantom{J(k,\xi;{u})=}+\sum_{\ell=t+d}^{N-1}\mathbb{E}\left\{\left[\begin{array}{cc}
\mathbb{E}_{\ell-d}X_\ell\\
u_\ell
\end{array}\right]^T\left[
\begin{array}{cc}
-P^{(d)}_\ell&H_\ell^T\\H_\ell&W_\ell
\end{array} \right] \left[\begin{array}{cc}
\mathbb{E}_{\ell-d}X_\ell\\
u_\ell
\end{array}\right]\right\}\nonumber \\
&&\hpm{J(k,\xi;{u})=}+\sum_{\ell=t+1}^{t+d-1}\mathbb{E}\left\{\left[\begin{array}{cc}
\mathbb{E}_{t}X_\ell\\
u_\ell
\end{array}\right]^T\left[
\begin{array}{cc}
A_\ell^TP^{(\ell+1-t)}_{\ell+1}A_\ell-P^{(\ell-t)}_\ell&H_\ell^T\\H_\ell&W_\ell
\end{array} \right]\left[\begin{array}{cc}
\mathbb{E}_{t}X_\ell\\
u_\ell
\end{array}\right]\right\}\nonumber \\
&&\hpm{J(k,\xi;{u})=}+\mathbb{E}\left\{\left[\begin{array}{cc}
X_t\\
u_t
\end{array}\right]^T\left[
\begin{array}{cc}
Q_t+A_t^T\big{(}P^{(0)}_{t+1}+P^{(1)}_{t+1}\big{)}A_t^T+C_t^TP^{(0)}_{t+1}C_t -P_{t}^{(0)} &H_t^T\\H_t&W_t
\end{array} \right]\left[\begin{array}{cc}
X_t\\
u_t
\end{array}\right]\right\}\nonumber \\
&&\hphantom{J(k,\xi;{u}) =}+\xi^TP^{(0)}_t\xi+\mathbb{E}\big{[}X_N^T\big{(}G-P_N^{(0)}\big{)}X_N \big{]}\\[1mm]
%
&&\hphantom{J(k,\xi;{u})}\geq \xi^TP^{(0)}_t\xi.
%
\end{eqnarray*}

\end{lemma}

\begin{remark}
The LMEIs (\ref{Riccati-inequality-1})-(\ref{Riccati-inequality-3}) are such constructed that the inequalities of Lemma \ref{Lemma-add-subtr} are satisfied. In this case, Problem (LQ) will be finite. Note that the LMEIs contain equality constraints; such new feature does not appear in deterministic LQ problems (with time delay) and standard stochastic LQ problems.

\end{remark}

Based on above preparations, we have the following theorem, which gives several equivalent characterizations on the solvability of Problem (LQ).

\begin{theorem}\label{Theorem-all-ii}
The following statements are equivalent.

(i) Problem (LQ) is finite.

(ii) Problem (LQ) is solvable.

(iii) The solution of (\ref{Riccati-1})-(\ref{Riccati-3}) has the property $W_kW_k^\dagger H_k=H_k, W_k\geq 0, k\in \mathbb{T}_t$, namely, the constrained Riccati-like equation set %
\begin{eqnarray}\label{Riccati-1-2}
\left\{
\begin{array}{l}
P^{(0)}_k=Q_k+A_k^T\big{(}P^{(0)}_{k+1}+P^{(1)}_{k+1}\big{)}A_k+C_k^TP^{(0)}_{k+1}C_k,\\[1mm]
P^{(i)}_k=A_k^TP^{(i+1)}_{k+1} A_k,~~i=1,...,d-1,\\[1mm]
P^{(d)}_k=-H_k^TW_k^\dagger H_k,\\[1mm]
P^{(0)}_N=G,~~P^{(j)}_N=0,~~j=1,...,d,\\[1mm]
W_kW_k^\dagger H_k=H_k, W_k\geq 0,\\[1mm]
k\in \mathbb{T}_{t+d},
\end{array}
\right.
\end{eqnarray}
\begin{eqnarray}\label{Riccati-2-2}
\left\{
\begin{array}{l}
P^{(0)}_k=Q_k+A_k^T\big{(}P^{(0)}_{k+1}+P^{(1)}_{k+1}\big{)}A_k+C_k^TP_{k+1}^{(0)}C_k,\\[1mm]
P^{(i)}_{k}=A_k^TP^{(i+1)}_{k+1}A_k,~~i=1,...,k-t-1,\\[1mm]
P^{(k-t)}_k=A_k^TP^{(k+1-t)}_{k+1}A_k-H_k^TW_k^\dagger H_k,\\[1mm]
W_kW_k^\dagger H_k=H_k, W_k\geq 0,\\[1mm]
k\in \{t+2,...,t+d-1\},
\end{array}
\right.
\end{eqnarray}
and
\begin{eqnarray}\label{Riccati-3-2}
\left\{
\begin{array}{l}
\left\{
\begin{array}{l}
P^{(0)}_{t+1}=Q_{t+1}+A_{t+1}^T\big{(}P^{(0)}_{t+2}+P^{(1)}_{t+2}\big{)}A_{t+1}+C_{t+1}^TP_{t+2}^{(0)}C_{t+1},\\[1mm]
%
%
P^{(1)}_{t+1}=A_{t+1}^TP^{(2)}_{{t+2}}A_{t+1}-H_{t+1}^TW_{t+1}^\dagger H_{t+1},
\end{array}
\right.\\[2mm]
P^{(0)}_t=Q_{t}+A_{t}^T\big{(}P^{(0)}_{t+1}+P^{(1)}_{t+1}\big{)}A_{t}+C_{t}^TP_{t+1}^{(0)}C_{t}-H_{t}^TW_{t}^\dagger H_{t},
\\[1mm]
W_kW_k^\dagger H_k=H_k, W_k\geq 0,~ k= t,t+1
\end{array}
\right.
\end{eqnarray}
are solvable in the sense
$$W_kW_k^\dagger H_k=H_k, W_k\geq 0, ~~k\in \mathbb{T}_t.$$

(iv) $\mathcal{M}$ defined in (\ref{M}) is nonempty.

Under any of the above conditions, an optimal control of Problem (LQ) for the initial pair $(k,\xi)$ is given by
\begin{eqnarray}\label{Theorem-all-u}
u_\ell^{k,\xi,*}=
-W_\ell^\dagger H_\ell\mathbb{E}_{\ell-d}X^{k,\xi,*}_\ell,~~\ell\in \mathbb{T}_{k}
\end{eqnarray}
with
\begin{eqnarray*}
\left\{
\begin{array}{l}
X^{k,\xi,*}_{\ell+1}=\big{(}A_{\ell}X^{k,\xi,*}_\ell-B_{\ell}W_\ell^\dagger H_\ell\mathbb{E}_{\ell-d}X^{k,\xi,*}_\ell\big{)}+\big{(}C_{\ell}X^{k,\xi,*}_\ell-D_{\ell}W_\ell^\dagger H_\ell\mathbb{E}_{\ell-d}X^{k,\xi,*}_\ell\big{)}w_\ell,\\[1mm]
X^{k,\xi,*}_k=\xi,~~\ell\in \mathbb{T}_{k}.
\end{array}
\right.
\end{eqnarray*}
Furthermore, the corresponding optimal value is
\begin{eqnarray}\label{Theorem-all-value}
V(k,\xi)=\left\{
\begin{array}{l}
\sum_{i=0}^{k-t}\xi^TP_k^{(i)}\xi,~~~k\in \{t,...,t+d-1\},\\[1mm]
\sum_{i=0}^{d}\xi^TP_k^{(i)}\xi,~~~k\in \mathbb{T}_{t+d}.
\end{array}
\right.
\end{eqnarray}

\end{theorem}

\emph{Proof}. See Appendix C. \endpf

\begin{corollary}\label{corollary-nonnegative}
Let $Q_k\geq 0, R_k\geq 0, G\geq 0, k\in \mathbb{T}_t$. Then, Problem (LQ) is solvable, and the corresponding Riccati-like equations (\ref{Riccati-1-2})-(\ref{Riccati-3-2}) are solvable.

\end{corollary}

\emph{Proof}. In this case, Problem (LQ) is finite, and the conclusion follows from Theorem  \ref{Theorem-all-ii}.   \endpf

\subsection{From the LMEIs to the Riccati-like equations\\------construct the solution of (\ref{Riccati-1-2})-(\ref{Riccati-3-2})}

In this subsection, a procedure is presented to construct the solution of Riccati-like equation set (\ref{Riccati-1-2})-(\ref{Riccati-3-2}) from an element $(\widetilde{P}^{(0)},...\widetilde{P}^{(d)})$ of $\mathcal{M}$ (i.e., a solution of the LMEIs (\ref{Riccati-inequality-1})-(\ref{Riccati-inequality-3})).  To do so, we introduce an auxiliary LQ problem.

Specifically, introduce the following weighting matrices
\begin{eqnarray}\label{weithing-matrices}
\left\{
\begin{array}{l}
\widetilde{Q}_k=Q_k+A_k^T\big{(}\widetilde{P}^{(0)}_{k+1}+\widetilde{P}^{(1)}_{k+1}\big{)}A_k+C_k^T\widetilde{P}^{(0)}_{k+1}C_k-\widetilde{P}^{(0)}_k,~~~ k\in \mathbb{T}_t,\\[1mm]
\widetilde{L}_k=\widetilde{H}_k^T\triangleq \ds\left\{
\begin{array}{ll}
\sum_{i=0}^{k+1-t}B_k^T \widetilde{P}^{(i)}_{k+1}A_k+D_k^T\widetilde{P}^{(0)}_{k+1}C_k,&~~k\in \{t,...,t+d-1\},\\[2mm]
\sum_{i=0}^{d}B_k^T \widetilde{P}^{(i)}_{k+1}A_k+D_k^T\widetilde{P}^{(0)}_{k+1}C_k,&~~k\in \mathbb{T}_{t+d},
\end{array}
\right.\\[1mm]
\widetilde{R}_k=\widetilde{W}_k\triangleq \left\{
\begin{array}{ll}
R_k+\sum_{i=0}^{k+1-t}B_k^T \widetilde{P}^{(i)}_{k+1} B_k+D_k^T\widetilde{P}^{(0)}_{k+1} D_k,&~~k\in \{t,...,t+d-1\},\\[2mm]
 R_k+\sum_{i=0}^{d}B_k^T \widetilde{P}^{(i)}_{k+1} B_k+D_k^T\widetilde{P}^{(0)}_{k+1} D_k,&~~k\in \mathbb{T}_{t+d},
\end{array}
\right.\\[1mm]
\widetilde{G}=G-\widetilde{P}^{(0)}_{N}.
\end{array}
\right.
\end{eqnarray}
Furthermore, for each $(k,\xi)\in \mathbb{T}_t\times \mathbb{R}^n$, let $X$ be the solution of (\ref{system-k}) and introduce the cost functional $\widetilde{J}(k,\xi; u)$ according to three different situations. Case 1: $k\in \mathbb{T}_{t+d}$, let
\begin{eqnarray}\label{cost-k---2}
&&\widetilde{J}(k,\xi; u)
 =\sum_{\ell=k}^{N-1}\mathbb{E}\big{[}X_\ell^T\widetilde{Q}_{\ell}X_\ell+2X_\ell^T\widetilde{L}_\ell u_\ell+u_{\ell}^T\widetilde{R}_{\ell}u_{\ell}\big{]}+\mathbb{E}\big{[}X_N^T\widetilde{G}X_N\big{]}\nonumber \\
 &&\hpm{\widetilde{J}(k,\xi; u)
 =}+\sum_{\ell=k}^{N-1}\mathbb{E}\big{[}-(\mathbb{E}_{\ell-d}X_\ell)^T\widetilde{P}^{(d)}_{\ell}\mathbb{E}_{\ell-d}X_\ell
\big{]};\end{eqnarray}
Case 2: $k\in \{t+1,...,t+d-1\}$, let
\begin{eqnarray}\label{cost-k---3}
&&\widetilde{J}(k,\xi; u)
 =\sum_{\ell=k}^{N-1}\mathbb{E}\big{[}X_\ell^T\widetilde{Q}_{\ell}X_\ell+2X_\ell^T\widetilde{L}^T_\ell u_\ell+u_{\ell}^T\widetilde{R}_{\ell}u_{\ell}\big{]}+\mathbb{E}\big{[}X_N^T\widetilde{G}X_N\big{]}\nonumber \\
&&\hphantom{\widetilde{J}(k,\xi; u)=} +\sum_{\ell=t+d}^{N-1}\mathbb{E}\big{[}-(\mathbb{E}_{\ell-d}X_\ell)^T\widetilde{P}^{(d)}_{\ell}\mathbb{E}_{\ell-d}X_\ell
\big{]}\nonumber \\
&&\hphantom{\widetilde{J}(k,\xi; u)=}+\sum_{\ell=k}^{t+d-1}\mathbb{E}\big{[}(\mathbb{E}_{\ell-d}X_\ell)^T\big{(}A_\ell^T\widetilde{P}^{(\ell+1-t)}_{\ell+1}A_\ell-\widetilde{P}^{(\ell-t)}_\ell\big{)}\mathbb{E}_{\ell-d}X_\ell \big{]};
\end{eqnarray}
Case 3: $k=t$, let
\begin{eqnarray}\label{cost-k---4}
&&\widetilde{J}(t,\xi; u)
 =\sum_{\ell=t}^{N-1}\mathbb{E}\big{[}X_\ell^T\widetilde{Q}_{\ell}X_\ell+2X_\ell^T\widetilde{L}^T_\ell u_\ell+u_{\ell}^T\widetilde{R}_{\ell}u_{\ell}\big{]}+\mathbb{E}\big{[}X_N^T\widetilde{G}X_N\big{]}\nonumber \\
&&\hphantom{\widetilde{J}(k,\xi; u)=} +\sum_{\ell=t+d}^{N-1}\mathbb{E}\big{[}-(\mathbb{E}_{\ell-d}X_\ell)^T\widetilde{P}^{(d)}_{\ell}\mathbb{E}_{\ell-d}X_\ell
\big{]}\nonumber \\
&&\hphantom{\widetilde{J}(k,\xi; u)=}+\sum_{\ell=t+1}^{t+d-1}\mathbb{E}\big{[}(\mathbb{E}_{\ell-d}X_\ell)^T\big{(}A_\ell^T\widetilde{P}^{(\ell+1-t)}_{\ell+1}A_\ell-\widetilde{P}^{(\ell-t)}_\ell\big{)}\mathbb{E}_{\ell-d}X_\ell \big{]}.
\end{eqnarray}
Corresponding to the above cost functional (\ref{cost-k---2})-(\ref{cost-k---4}), the system (\ref{system-k}) and the admissible control set (\ref{admissible-control--k-1})-(\ref{admissible-control--k-2}), we denote such an LQ problem as Problem (LQ)$_a$ for the initial pair $(k,\xi)$.

The cost functional $\widetilde{J}(k,\xi;u)$ is such constructed in (\ref{cost-k---2})-(\ref{cost-k---4}) that it is finite from below. This is proved in the following proposition.

\begin{proposition}\label{Lemma-axuxiliary}
For any $(k,\xi)\in \mathbb{T}_t\times \mathbb{R}^n$, $\widetilde{J}(t,\xi; u)\geq 0$. Hence, Problem (LQ)$_a$ is finite.

\end{proposition}

\emph{Proof}. For (\ref{cost-k---4}), we have
\begin{eqnarray*}\label{finite-f-3}
&&\hspace{-2em}\widetilde{J}(t,\xi; u)
 =\sum_{\ell=t+d}^{N-1}\mathbb{E}\left\{X_\ell^T\widetilde{Q}_{\ell}X_\ell+\left[
\begin{array}{c}
\mathbb{E}_{\ell-d}X_\ell\\
u_\ell
\end{array}
\right]^T\left[
\begin{array}{cc}
-\widetilde{P}^{(d)}_\ell&\widetilde{H}_\ell^T\\\widetilde{H}_\ell&\widetilde{W}_\ell
\end{array} \right]\left[
\begin{array}{c}
\mathbb{E}_{\ell-d}X_\ell\\
u_\ell
\end{array}
\right]\right\}\nonumber \\
 &&\hspace{-2em}\hphantom{\widetilde{J}(k,\xi; u)
 =}+\sum_{\ell=t+1}^{t+d-1}\mathbb{E}\left\{X_\ell^T\widetilde{Q}_{\ell}X_\ell+\left[
\begin{array}{c}
\mathbb{E}_{\ell-d}X_\ell\\
u_\ell
\end{array}
\right]^T\left[
\begin{array}{cc}
A_\ell^T\widetilde{P}^{(\ell+1-t)}_{\ell+1}A_\ell-\widetilde{P}^{(\ell-t)}_\ell&\widetilde{H}_\ell^T\\\widetilde{H}_\ell&\widetilde{W}_\ell
\end{array} \right]\left[
\begin{array}{c}
\mathbb{E}_{\ell-d}X_\ell\\
u_\ell
\end{array}
\right]\right\}   \nonumber \\
&&\hspace{-2em}\hphantom{\widetilde{J}(k,\xi; u)
 =}+\mathbb{E}\left\{
\left[
\begin{array}{c}
X_t\\u_t
\end{array}
\right]^T\left[
\begin{array}{cc}
\widetilde{Q}_t&\widetilde{H}_t^T\\\widetilde{H}_t&\widetilde{W}_t\end{array}
\right]   \left[
\begin{array}{c}
X_t\\u_t
\end{array}
\right]
 \right\}+\mathbb{E}\big{[}X_N^T\widetilde{G}X_N\big{]}\nonumber\\[1mm]
&&\hspace{-2em}\hphantom{\widetilde{J}(k,\xi; u)} \geq 0.
\end{eqnarray*}
The inequality above is due to the fact $(\widetilde{P}^{(0)},...,\widetilde{P}^{(d)})\in \mathcal{M}$. Similarly, we can prove other cases. Hence, $\widetilde{J}(k,\xi; u)\geq 0$ for any $(k,\xi)\in \mathbb{T}_t\times \mathbb{R}^n$.  \endpf

Let us make some observations about $\widetilde{J}(t,\xi;u)$. By adding to and subtracting
\begin{eqnarray*}
&&\sum_{k=t}^{N-1}\mathbb{E}\Big{\{}\sum_{i=0}^{d}(\mathbb{E}_{k+1-i}X_{k+1})^TU^{(i)}_{k+1}\mathbb{E}_{k+1-i}X_{k+1}-\sum_{i=0}^{d}(\mathbb{E}_{k-i}X_{k})^TU^{(i)}_{k}\mathbb{E}_{k-i}X_{k}\Big{\}}\\
&&+\sum_{k=t}^{t+d-1}\mathbb{E}\Big{\{}\sum_{i=0}^{k+1-t}(\mathbb{E}_{k+1-i}X_{k+1})^TU^{(i)}_{k+1}\mathbb{E}_{k+1-i}X_{k+1}-\sum_{i=0}^{k-t}(\mathbb{E}_{k-i}X_{k})^TU^{(i)}_{k}\mathbb{E}_{k-i}X_{k}\Big{\}}
\end{eqnarray*}
from $\widetilde{J}(t,\xi;u)$, we have
\begin{eqnarray}\label{add-subtr---new}
&&\hspace{-2em}\widetilde{J}(t,\xi;{u}) =\sum_{k=t+d}^{N-1}\mathbb{E}\Big{\{}X_k^T\big{[}\widetilde{Q}_k+A_k^T\big{(}U^{(0)}_{k+1}+U^{(1)}_{k+1}\big{)}A_k^T+C_k^TU^{(0)}_{k+1}C_k -U_{k}^{(0)}\big{]}X_k \nonumber  \\
&&\hspace{-2em}\hphantom{J(t,\xi;{u}) =}+\sum_{i=1}^{d-1}(\mathbb{E}_{k-i}X_k)^T\big{[}A_k^TU^{(i+1)}_{k+1}A_k-U_k^{(i)} \big{]}\mathbb{E}_{k-i}X_k-(\mathbb{E}_{k-d}X_k)^T\big{(}U^{(d)}_k+\widetilde{P}_{k}^{(d)}\big{)}\mathbb{E}_{k-d}X_k\nonumber  \\
&&\hspace{-2em}\hphantom{J(t,\xi;{u}) =}+2(\mathcal{H}_k\mathbb{E}_{k-d}X_k)^Tu_k+u_k^T\mathcal{W}_ku_k \Big{\}}\nonumber  \\
&&\hspace{-2em}\hphantom{J(t,\xi;{u}) =}+\sum_{k=t+2}^{t+d-1}\mathbb{E}\Big{\{}X_k^T\big{[}\widetilde{Q}_k+A_k^T\big{(}U^{(0)}_{k+1}+U^{(1)}_{k+1}\big{)}A_k^T+C_k^TU^{(0)}_{k+1}C_k -U_{k}^{(0)}\big{]}X_k \nonumber  \\
&&\hspace{-2em}\hphantom{J(t,\xi;{u}) =}+\sum_{i=1}^{k-t-1}(\mathbb{E}_{k-i}X_k)^T\big{[}A_k^TU^{(i+1)}_{k+1}A_k-U_k^{(i)} \big{]}\mathbb{E}_{k-i}X_k\nonumber\\
&&\hspace{-2em}\hphantom{J(t,\xi;{u}) =}+(\mathbb{E}_{t}X_k)^T\big{[}A_k^T\big{(}U^{(k+1-t)}_{k+1}+\widetilde{P}^{(k+1-t)}_{k+1}\big{)}A_k-\big{(}U^{(k-t)}_k+\widetilde{P}_k^{(k-t)}\big{)}\big{]}\mathbb{E}_{t}X_k\nonumber\\
&&\hspace{-2em}\hphantom{J(t,\xi;{u}) =}+2(\mathcal{H}_k\mathbb{E}_{t}X_k)^Tu_k+u_k^T\mathcal{W}_ku_k \Big{\}}\nonumber  \\[1mm]
&&\hspace{-2em}\hphantom{J(t,\xi;{u}) =}+\mathbb{E}\Big{\{}X_{t+1}^T\big{[}\widetilde{Q}_{t+1}+A_{t+1}^T\big{(}U^{(0)}_{{t+2}}+U^{(1)}_{{t+2}}\big{)}A_{t+1}^T+C_{t+1}^TU^{(0)}_{{t+2}}C_{t+1} -U_{{t+1}}^{(0)}\big{]}X_{t+1} \nonumber  \\[1mm]
&&\hspace{-2em}\hphantom{J(t,\xi;{u}) =}+(\mathbb{E}_{t}X_{t+1})^T\big{[}A_{t+1}^T\big{(}U^{(2)}_{t+2}+\widetilde{P}^{(2)}_{t+2}\big{)}A_{t+1}-\big{(}U^{(1)}_{t+1}+\widetilde{P}_{t+1}^{(1)}\big{)}\big{]} \mathbb{E}_{t}X^0_{t+1}+2(\mathcal{H}_{t+1}\mathbb{E}_{t}X_{t+1})^Tu_{t+1}\nonumber  \\[1mm]
&&\hspace{-2em}\hphantom{J(t,\xi;{u}) =}+u_{t+1}^T\mathcal{W}_{t+1}u_{t+1} \Big{\}}+\mathbb{E}\Big{\{}X_{t}^T\big{[}\widetilde{Q}_{t}+A_{t}^T\big{(}U^{(0)}_{{t+1}}+U^{(1)}_{{t+1}}\big{)}A_{t}^T+C_{t}^TU^{(0)}_{{t+1}}C_{t} -U_{{t}}^{(0)}\big{]}X_{t} \nonumber  \\[1mm]
%
%
%
&&\hspace{-2em}\hphantom{J(t,\xi;{u}) =}+2(\mathcal{H}_{t}X_{t})^Tu_{t}+u_{t}^T\mathcal{W}_{t}u_{t} \Big{\}}+\xi^TU^{(0)}_t\xi.
%
%
%
\end{eqnarray}
In the above, $(U^{(0)},...,U^{(d)})$ is to be determined and
\begin{eqnarray}\label{W-mathcal}
\mathcal{W}_{k}=\ds\left\{
\begin{array}{ll}
\widetilde{R}_k+\sum_{i=0}^{k+1-t}B_k^T U^{(i)}_{k+1} B_k+D_k^TU^{(0)}_{k+1} D_k,&~~k\in \{t,...,t+d-1\},\\[2mm]
\widetilde{R}_k+\sum_{i=0}^{d}B_k^T U^{(i)}_{k+1} B_k+D_k^TU^{(0)}_{k+1} D_k,&~~k\in \mathbb{T}_{t+d},
\end{array}
\right.
\end{eqnarray}
and
\begin{eqnarray}\label{H-mathcal}
\mathcal{H}_{k}=\ds\left\{
\begin{array}{ll}
\sum_{i=0}^{k+1-t}B_k^T U^{(i)}_{k+1}A_k+D_k^TU^{(0)}_{k+1}C_k+\widetilde{L}_k^T,&~~k\in \{t,...,t+d-1\},\\[2mm]
\sum_{i=0}^{d}B_k^T U^{(i)}_{k+1}A_k+D_k^TU^{(0)}_{k+1}C_k+\widetilde{L}_k^T,&~~k\in \mathbb{T}_{t+d}.
\end{array}
\right.
\end{eqnarray}
In fact, introduce the Riccati-like equation set
\begin{eqnarray}\label{Riccati-1-3}
\left\{
\begin{array}{l}
U^{(0)}_k=\widetilde{Q}_k+A_k^T\big{(}U^{(0)}_{k+1}+U^{(1)}_{k+1}\big{)}A_k+C_k^TU^{(0)}_{k+1}C_k,\\[1mm]
U^{(i)}_k=A_k^TU^{(i+1)}_{k+1} A_k,~~i=1,...,d-1,\\[1mm]
U^{(d)}_k=-\widetilde{P}^{(d)}_k-\mathcal{H}_k^T\mathcal{W}_k^\dagger \mathcal{H}_k,\\[1mm]
U^{(0)}_N=\widetilde{G},~~U^{(j)}_N=0,~~j=1,...,d,\\[1mm]
k\in \mathbb{T}_{t+d},
\end{array}
\right.
\end{eqnarray}
\begin{eqnarray}\label{Riccati-2-3}
\left\{
\begin{array}{l}
U^{(0)}_k=\widetilde{Q}_k+A_k^T\big{(}U^{(0)}_{k+1}+U^{(1)}_{k+1}\big{)}A_k+C_k^TU_{k+1}^{(0)}C_k,\\[1mm]
U^{(i)}_{k}=A_k^TU^{(i+1)}_{k+1}A_k,~~i=1,...,k-t-1,\\[1mm]
U^{(k-t)}_k=A_k^T(\widetilde{P}^{(k+1-t)}_{k+1}+U^{(k+1-t)}_{k+1})A_k-\widetilde{P}^{(k-t)}_k-\mathcal{H}_k^T\mathcal{W}_k^\dagger \mathcal{H}_k,\\[1mm]
k\in \{t+2,...,t+d-1\},
\end{array}
\right.
\end{eqnarray}
and
\begin{eqnarray}\label{Riccati-3-3}
\left\{
\begin{array}{l}
\left\{
\begin{array}{l}
U^{(0)}_{t+1}=\widetilde{Q}_{t+1}+A_{t+1}^T\big{(}U^{(0)}_{t+2}+U^{(1)}_{t+2}\big{)}A_{t+1}+C_{t+1}^TU_{t+2}^{(0)}C_{t+1},\\[1mm]
%
%
U^{(1)}_{t+1}=A_{t+1}^T(\widetilde{P}^{(2)}_{t+2}+U^{(2)}_{t+2})A_{t+1}-\widetilde{P}^{(1)}_{t+1}-\mathcal{H}_{t+1}^T\mathcal{W}_{t+1}^\dagger \mathcal{H}_{t+1},\\[1mm]
\end{array}
\right.\\
U^{(0)}_{t}=\widetilde{Q}_{t}+A_{t}^T\big{(}U^{(0)}_{t+1}+U^{(1)}_{t+1}\big{)}A_{t}+C_{t}^TU_{t+1}^{(0)}C_{t}-\mathcal{H}_{t}^T\mathcal{W}_{t}^\dagger \mathcal{H}_{t},\\[1mm]
\end{array}
\right.
\end{eqnarray}
with $\mathcal{W}_k, \mathcal{H}_k$ being given by (\ref{W-mathcal}) (\ref{H-mathcal}); by analysis similar to (\ref{add-subtr---new}), we then have the following result.

\begin{lemma}\label{Lemma-add-subtr-2}
Let $(U^{(0)},...,U^{(d)})$ be the solution of (\ref{Riccati-1-3})-(\ref{Riccati-3-3}). Then,
\begin{eqnarray*}
\widetilde{J}(k,\xi;u)=\sum_{\ell=k}^{N-1}\mathbb{E}\Big{\{}(\mathbb{E}_{\ell-d}X_\ell)^T\mathcal{H}_\ell^T\mathcal{W}^\dagger_\ell \mathcal{H}_\ell\mathbb{E}_{\ell-d}X_\ell+2(\mathcal{H}_\ell\mathbb{E}_{\ell-d}X_\ell)^Tu_k+u_k^T\mathcal{W}_ku_k \Big{\}}+\widetilde{\Pi}_k(\xi),
\end{eqnarray*}
where
\begin{eqnarray*}
\widetilde{\Pi}_k(\xi)=\left\{
\begin{array}{l}
\sum_{i=0}^{k-t}\xi^TU_k^{(i)}\xi,~~~k\in \{t,...,t+d-1\},\\
\sum_{i=0}^{d}\xi^TU_k^{(i)}\xi,~~~k\in \mathbb{T}_{t+d}.
\end{array}
\right.
\end{eqnarray*}

\end{lemma}

Based on what we have prepared above, we can construct a solution of (\ref{Riccati-1-2})-(\ref{Riccati-3-2}) from $(\widetilde{P}^{(0)},...,\widetilde{P}^{(d)})\in \mathcal{M}$.

\begin{theorem}\label{Theorem--construction}
The following statements hold.

(i) The solution of Riccati-like equation set (\ref{Riccati-1-3})-(\ref{Riccati-3-3}) has the property
\begin{eqnarray*}
\mathcal{W}_k\geq 0,~~\mathcal{W}_k\mathcal{W}^\dagger_k\mathcal{H}_k=\mathcal{H}_k,~~k\in \mathbb{T}_t.
\end{eqnarray*}

(ii) Let $P_k^{(i)}=\widetilde{P}_k^{(i)}+U_k^{(i)}, k\in \mathbb{T}_t, i=0,...,d$. Then, such a $(P^{(0)},...,P^{(d)})$ is a solution of the constrained Riccati-like equation set (\ref{Riccati-1-2})-(\ref{Riccati-3-2}).

\end{theorem}

\emph{Proof}.
From Proposition \ref{Lemma-axuxiliary}, Problem (LQ)$_a$ is finite for any initial pair $(k,\xi)\in \mathbb{T}_t\times \mathbb{R}^n$; hence it is solvable.
Combining Lemma \ref{Lemma-add-subtr-2} and the part of proving the equivalence between (i) and (iii) of Theorem \ref{Theorem-all-ii}, we must have (i) of this theorem.
(ii) follows from some simple calculations. \endpf

\begin{remark}
By Theorem \ref{Theorem--construction}, we can construct a solution of the constrained Riccati-like equation set from a solution of the LMEIs.
This result is  potentially useful to study the algebraic Riccati-like equations that we will encounter in the infinite-horizon version of Problems (LQ). For more about standard infinite-horizon stochastic LQ problems, we can refer to, for example, \cite{Ait-Moore-Zhou} \cite{Yao-Zhang-Zhou}.

\end{remark}

\subsection{The unique solvability of Problem (LQ)}

In the following, we will study the uniform convexity of the cost functional, which is motivated by some results of \cite{Sun-Yong-Siam-2016}.
The functional $u\mapsto J(t,x;u)$ is called uniformly convex if there exists a $\lambda>0$ such that for any $u\in \mathcal{U}_{ad}^t$
\begin{eqnarray}\label{strict-convex}
J(t,0;u)\geq \lambda ||u||^2=\lambda \sum_{k=t}^{N-1}\mathbb{E}|u_k|^2.
\end{eqnarray}
From Proposition \ref{proposition-operator-Hilbert}, Problem (LQ) will have a unique optimal control if $u\mapsto J(t, 0;u)$ is uniformly convex.
%


\begin{lemma}\label{Lemma-norm}
For $\Phi\in L^2(\mathbb{T}_t; \mathbb{R}^{m\times n})$ and (\ref{system-0}) with $u\in \mathcal{U}_{ad}^t$, there exist $\gamma_1, \gamma_2$ with property $0<\gamma_2<\gamma_1$ such that
\begin{eqnarray}\label{Lemma-norm-1}
\gamma_2 \sum_{k=t}^{N-1}\mathbb{E}|u_k|^2\leq \sum_{k=t}^{N-1}\mathbb{E}|u_k-\Phi_k \mathbb{E}_{k-d}X^0_k|^2\leq \gamma_1 \sum_{k=t}^{N-1}\mathbb{E}|u_k|^2.
\end{eqnarray}
%

\end{lemma}

\emph{Proof}. For $\Phi\in L^2(\mathbb{T}_t; \mathbb{R}^{m\times n})$, define a bounded linear operator for $\mathcal{U}_{ad}^t$ to $\mathcal{U}_{ad}^t$
\begin{eqnarray}\label{operator-O}
\mathcal{O}u=u-\Phi \mathbb{E}_{\cdot-d}X^{0},
\end{eqnarray}
where ${u}-\Phi\mathbb{E}_{\cdot-d}X^0$ is the control $\{{u}_k-\Phi_k\mathbb{E}_{k-d}X^0_k, ~~k\in \mathbb{T}_t\}$.
Note that $\mathcal{O}u=0$ implies $u=0$, i.e., $\mathcal{O}$ is an injection. Let
\begin{eqnarray*}
p_\Phi(u)=||\mathcal{O}u||=\sqrt{\sum_{k=t}^{N-1}\mathbb{E}|u_k-\Phi_k \mathbb{E}_{k-d}X^0_k|^2},
\end{eqnarray*}
which is indeed a norm on $\mathcal{U}_{ad}^t$. Furthermore, for any given $u^{(n)}\in \mathcal{U}^t_{ad}$, we have when $n\mapsto \infty$
\begin{eqnarray*}
p_{\Phi}(u^{(n)})\mapsto 0~~ \Leftrightarrow~~||u^{(n)}||=\sqrt{\sum_{k=t}^{N-1} \mathbb{E}|u_k^{(n)}|^2}\mapsto 0.
\end{eqnarray*}
Therefore, $p_{\Phi}(\,\cdot\,)$ is equivalent to the norm $||\cdot||$ on  $\mathcal{U}_{ad}^t$. We then claim (\ref{Lemma-norm-1}).\endpf 

\begin{theorem}\label{Theorem-nec-suff}
The following statements are equivalent.

(i) Problem (LQ) is uniquely solvable at the initial pair $(t,x)$.

(ii) Riccati-like equation set (\ref{Riccati-1-2})-(\ref{Riccati-3-2}) is solvable, and $W_k>0, k\in \mathbb{T}_t$.

(iii) $u \mapsto J(t,x; u)$ is uniformly convex for $u\in \mathcal{U}^t_{ad}$.

(iv) For any $k\in \mathbb{T}_t$, $u \mapsto J(k,\xi; u)$ is uniformly convex for $u\in \mathcal{U}^k_{ad}$.

(v) Problem (LQ) is uniquely solvable at any initial pair $(k,\xi)\in \mathbb{T}_t\times \mathbb{R}^n$.

Under any of the above conditions, the optimal control of Problem (LQ) for the initial pair $(k,\xi)$ is given by
\begin{eqnarray}\label{Theorem-nec-suff-u}
u_\ell^{k,\xi,*}=
-W_\ell^{-1} H_\ell\mathbb{E}_{\ell-d}X^{k,\xi,*}_\ell,~~\ell\in \mathbb{T}_{k}
\end{eqnarray}
with $X^{k,\xi,*}$ given by
\begin{eqnarray*}
\left\{
\begin{array}{l}
X^{k,\xi,*}_{\ell+1}=\big{(}A_{\ell}X^{k,\xi,*}_\ell-B_{\ell}W_\ell^{-1} H_\ell\mathbb{E}_{\ell-d}X^{k,\xi,*}_\ell\big{)}+\big{(}C_{\ell}X^{k,\xi,*}_\ell-D_{\ell}W_\ell^{-1} H_\ell\mathbb{E}_{\ell-d}X^{k,\xi,*}_\ell\big{)}w_\ell,\\[1mm]
X^{k,\xi,*}_k=\xi,~~\ell\in \mathbb{T}_{k}.
\end{array}
\right.
\end{eqnarray*}
%

\end{theorem}

\emph{Proof}. $(i)\Rightarrow (ii)$. This can be achieved by undating the proof of Theorem \ref{Theorem-necessary}. Let $u^{t,x,*}$ be the unique optimal control of Problem (LQ) for the initial pair $(t,x)$. Noting (\ref{Theorem-necessary-proof-3}) and that Theorem \ref{Theorem--Nece-Suff-fixed} is of necessary and sufficient conditions, we have
\begin{eqnarray*}\label{}
0=W_ku^{t,x,*}_k+H_kX_k^{t,x,*},~~k\in \mathbb{T}_t.
\end{eqnarray*}
As the optimal control uniquely exists, we must have that $W_k$ is nonsingular, $k\in \mathbb{T}_t$. Otherwise, any controls of the following form
\begin{eqnarray}\label{Theorem-nec-suff-1}
\widehat{{u}}^{t,x,*}_k=-W_k^\dagger H_k{\widehat{X}}_k^{t,x,*}+\big{(}I-W_k^\dagger W_k \big{)}\Upsilon_k,~~\Upsilon_k\in \mathcal{F}_{k-d},~~k\in \mathbb{T}_t
\end{eqnarray}
is also an optimal control, where
\begin{eqnarray*}
\left\{\begin{array}{l}
\widehat{X}^{t,x,*}_{k+1}=\big{(}A_{k}\widehat{X}^{t,x,*}_k+B_{k}\widehat{u}^{t,x,*}_{k}\big{)}+\big{(}C_{k}\widehat{X}^{t,x,*}_k+D_{k}\widehat{u}^{t,x,*}_{k}\big{)}w_k, \\[1mm]
\widehat{X}^{t,x,*}_t=x,~~k\in \mathbb{T}_t.
\end{array}
\right.
\end{eqnarray*}
%
Since $W_k, k\in \mathbb{T}_t$, are all invertible, from (\ref{add-subtr}) and (\ref{convexity-condition}) we have
\begin{eqnarray}\label{Theorem-nec-suff-2}
&&J(t,0;u)
%
%
%
%
%
=\sum_{k=t}^{N-1}(u_k+W^{-1}_kH_k\mathbb{E}_{k-d}X^0_k)^TW_k(u_k+W_k^{-1}H_k\mathbb{E}_{k-d}X^0_k)
\geq 0,
\end{eqnarray}
where $X^0$ is given in (\ref{system-0}).
%
Letting $\Phi_k=-W^{-1}_kH_k$, from Lemma \ref{Lemma-cost} we know that the linear operator $\mathcal{O}$ defined in (\ref{operator-O}) is a surjection from $\mathcal{U}^t_{ad}$ to $\mathcal{U}^t_{ad}$. 
%
%
%
%
%
%
Noting that (\ref{Theorem-nec-suff-2}) holds for any $u\in \mathcal{U}^t_{ad}$, we must have
\begin{eqnarray*}
W_k>0,~~k\in \mathbb{T}_t,
\end{eqnarray*}
which implies (ii).

$(ii)\Rightarrow (iii)$. Similarly to (\ref{Theorem-nec-suff-2}), it holds that
\begin{eqnarray*}
J(t,0;u)
=\sum_{k=t}^{N-1}(u_k+W^{-1}_kH_k\mathbb{E}_{k-d}X^0_k)^TW_k(u_k+W_k^{-1}H_k\mathbb{E}_{k-d}X^0_k)\nonumber.
%
%
%
\end{eqnarray*}
From Lemma \ref{Lemma-norm}, we have
\begin{eqnarray*}
J(t,0;u)\geq \lambda_{min}\sum_{k=t}^{N-1}\mathbb{E}|u_k+W^{-1}_kH_k\mathbb{E}_{t}X^0_k|^2 \geq  \lambda_{min}\gamma_2 \sum_{k=t}^{N-1}\mathbb{E}|u_k|^2,
\end{eqnarray*}
where $\lambda_{min}>0$ denotes the minimal eigenvalue among all the eigenvalues of $W_k, k\in \mathbb{T}_t$. Hence, $u \mapsto J(t,x; u)$ is uniformly convex.

(iii)$\Rightarrow$(iv). Let $u \mapsto J(t,x; u)$ be uniformly convex for $u\in \mathcal{U}^t_{ad}$. Now for any $u=(u_k,...,u_{N-1})\in \mathcal{U}_{ad}^k$, let $v=(0,...,0, u_k,...,u_{N-1})\in \mathcal{U}_{ad}^t$. Then, we have
\begin{eqnarray*}
J(k,0;u)=J(t,0;v)\geq \lambda \sum_{\ell=t}^{N-1}\mathbb{E}|v_\ell|^2= \lambda \sum_{\ell=k}^{N-1}\mathbb{E}|u_\ell|^2.
\end{eqnarray*}
Hence, $u \mapsto J(k,\xi; u)$ is uniformly convex.

(iv)$\Rightarrow$(v). From Proposition \ref{proposition-operator-Hilbert}, Problem (LQ) for the initial pair $(k,\xi)$ admits a unique optimal control.

(v)$\Rightarrow$(i). This is clear.

Under any of the above conditions, we have (\ref{Theorem-nec-suff-u}).  \endpf

\begin{remark}
The theorem above shows that Problem (LQ) is uniquely solvable at the initial pair $(t,x)$ if and only if Problem (LQ) is uniquely solvable at any initial pair $(k,\xi)\in \mathbb{T}_t\times \mathbb{R}^n$. This result links Section \ref{section--fixed} with this section.
Note, here, that the condition of uniform convexity plays a key role.
\end{remark}

\section{Example}\label{Section-example}

In this section, we shall present an example to illustrate the theory derived above.

\begin{example}
Consider a version of Problem (LQ) whose system matrices and weighting matrices are
\begin{eqnarray*}
&&\hspace{-2em}A_{0}=\left[
\begin{array}{cc}
-1.2 & 0.41\\
-0.3& 0.89
\end{array}
\right],
A_{1}=\left[
\begin{array}{cc}
2.32 & -0.35\\
0.31& 0.3
\end{array}
\right], A_{2}=\left[
\begin{array}{cc}
2.15 & -0.3\\
1.2& 4
\end{array}
\right], A_{3}=\left[
\begin{array}{cc}
-1.15 & -0.23\\
-2& 1
\end{array}
\right],\\[1mm]
&&\hspace{-2em}B_{0}=\left[
\begin{array}{cc}
2.25 & 0.6\\
-1.2& 3
\end{array}
\right], B_{1}=\left[
\begin{array}{cc}
2.2 & -1.32\\
0.5&  3
\end{array}
\right], B_{2}=\left[
\begin{array}{cc}
5.15 & 0\\
0& 5.6
\end{array}
\right], B_{3}=\left[
\begin{array}{cc}
1.35 & 1\\
-0.2& 1
\end{array}
\right],\\[1mm]
&&\hspace{-2em}C_{0}=\left[
\begin{array}{cc}
2.6 & 1\\
-1.73&  7.8
\end{array}
\right], C_{1}=\left[
\begin{array}{cc}
2.5 & 0.73\\
-1.47& 5.2
\end{array}
\right], C_{2}=\left[
\begin{array}{cc}
2.6 & 1.63\\
-1& 3.7
\end{array}
\right], C_{3}=\left[
\begin{array}{cc}
1.6 & 0.6\\
1& 2.1
\end{array}
\right],\\[1mm]
&&\hspace{-2em}
D_{0}=\left[
\begin{array}{cc}
2.4 & 1.93\\
1.07& 3
\end{array}
\right], D_{1}=\left[
\begin{array}{cc}
2.8 & 1.03\\
-1.23& 6
\end{array}
\right],
D_{2}=\left[
\begin{array}{cc}
0.5 & 0.2\\
1.1 & 2.65
\end{array}
\right], D_{3}=\left[
\begin{array}{cc}
1.5 & -1\\
-0.16 & 1.65
\end{array}
\right],\\[1mm]
&&\hspace{-2em}Q_{0}=\left[
\begin{array}{cc}
-2 & 0.8\\
0.8&  -1.6
\end{array}
\right], Q_{1}=\left[
\begin{array}{cc}
4 & 0\\
0& 0
\end{array}
\right], Q_{2}=\left[
\begin{array}{cc}
-0.5 & 0\\
0&  1
\end{array}
\right], Q_{3}=\left[
\begin{array}{cc}
1 & 0\\
0&  4
\end{array}
\right],\\[1mm]
&&\hspace{-2em}R_{0}=\left[
\begin{array}{cc}
-5 & 0\\
0&  -4
\end{array}
\right], R_{1}=\left[
\begin{array}{cc}
-2 & 0.1\\
0.1& 5
\end{array}
\right], R_{2}=\left[
\begin{array}{cc}
4 & -0.3\\
-0.3& 7
\end{array}
\right], R_{3}=\left[
\begin{array}{cc}
2 & -0.3\\
-0.3& 0
\end{array}
\right],\\
&&\hspace{-2em}G=\left[
\begin{array}{cc}
2 & -0.3\\
-0.3&  0
\end{array}
\right].
\end{eqnarray*}
Let $N=4$ and $d=2$ in (\ref{system-1}) in (\ref{admissible-control--t}). Find the optimal control.
\end{example}

In this case, the constrained Riccati-like equation set (\ref{Riccati-1-2})-(\ref{Riccati-3-2}) becomes to
\begin{eqnarray}\label{Riccati-1-4}
\left\{
\begin{array}{l}
P^{(0)}_k=Q_k+A_k^T\big{(}P^{(0)}_{k+1}+P^{(1)}_{k+1}\big{)}A_k+C_k^TP^{(0)}_{k+1}C_k,\\[1mm]
P^{(1)}_k=A_k^TP^{(i+1)}_{k+1} A_k,\\[1mm]
P^{(2)}_k=-H_k^TW_k^\dagger H_k,\\[1mm]
P^{(0)}_4=G,~~P^{(1)}_4=P^{(2)}_4=0, \\[1mm]
W_kW_k^\dagger H_k=H_k, W_k\geq 0,\\[1mm]
k\in \{2,3\},
\end{array}
\right.
\end{eqnarray}
and
\begin{eqnarray}\label{Riccati-3-4}
\left\{
\begin{array}{l}
\left\{
\begin{array}{l}
P^{(0)}_{1}=Q_{1}+A_{1}^T\big{(}P^{(0)}_{2}+P^{(1)}_{2}\big{)}A_{1}+C_{1}^TP_{2}^{(0)}C_{1},\\[1mm]
%
%
P^{(1)}_{1}=A_{1}^TP^{(2)}_{{2}}A_{1}-H_{1}^TW_{1}^\dagger H_{1},
\end{array}
\right.\\[2mm]
P^{(0)}_0=Q_{0}+A_{0}^T\big{(}P^{(0)}_{1}+P^{(1)}_{1}\big{)}A_{0}+C_{0}^TP_{1}^{(0)}C_{0}-H_{0}^TW_{0}^\dagger H_{0},
\\[1mm]
W_kW_k^\dagger H_k=H_k, W_k\geq 0,~ k= 0,1.
\end{array}
\right.
\end{eqnarray}
By some calculations, we have
\begin{eqnarray*}
&&W_0=\left[
\begin{array}{cc}
    7926 &   4307\\
    4307  &  1403
\end{array}
\right], W_1=\left[
\begin{array}{cc}
    749.8 &  -120.6\\
   -120.6  &  6637
\end{array}
\right],\\
&&W_2=\left[
\begin{array}{cc}
   28.8150 &   5.7102\\
    5.7102 & 151.0654
    \end{array}
\right], W_3=\left[
\begin{array}{cc}
   10.4510 &  -1.7355\\
   -1.7355 &   4.3900
    \end{array}
\right],
\end{eqnarray*}
which are positive definite. Hence, (\ref{Riccati-1-4})-(\ref{Riccati-3-4}) are solvable. Furthermore, the unique optimal control is given by
\begin{eqnarray*}
u_k^{0,x,*}=
-W_k^{-1} H_k\mathbb{E}_{k-2}X^{0,x,*}_k,~~k\in \{0,1,2,3\},
\end{eqnarray*}
where $-W_k^\dagger H_k, k=0,1,2,3$, are
\begin{eqnarray*}
&& -W_0^{-1} H_0=\left[
\begin{array}{cc}
   -1.5730 &   1.2102\\
    1.0877 &  -2.9347
\end{array}
\right],~~-W_1^{-1} H_1=\left[
\begin{array}{cc}
   -0.9460 &   0.0731\\
    0.0572 &  -0.8292
\end{array}
\right], \\
&&
-W_2^{-1} H_2=\left[
\begin{array}{cc}
   -0.3940 &  -0.5321\\
   -0.1330 &  -0.8525
    \end{array}
\right],~~
-W_3^{-1} H_3=\left[
\begin{array}{cc}
   -0.0069 &   0.0791\\
    1.1469 &   0.3861
    \end{array}
\right],
\end{eqnarray*}
and $X^{0,x,*}$ is given by
\begin{eqnarray*}\label{example--state}
\left\{
\begin{array}{l}
X^{0,x,*}_{k+1}=\big{(}A_{k}X^{0,x,*}_k-B_{k}W_k^{-1} H_k\mathbb{E}_{k-2}X^{0,x,*}_k\big{)}+\big{(}C_{k}X^{0,x,*}_k-D_{k}W_k^{-1} H_k\mathbb{E}_{k-2}X^{0,x,*}_k\big{)}w_k,\\[1mm]
X^{0,x,*}_0=x,~~k\in \{0,1,2,3\}.
\end{array}
\right.
\end{eqnarray*}

\section{Conclusion}\label{section-conclusion}

In this paper, an indefinite stochastic LQ problem with transmission delay and multiplicative noises is studied. Based on some abstract consideration, necessary and sufficient conditions are given, respectively, for the case with a fixed initial pair and the case with all the initial pairs. Further, a set of constrained discrete-time Riccati-like equations and a set of linear matrix equality-inequalities are introduced, which are used to characterize the existence of the delayed optimal control of Problem (LQ). Moreover, the unique solvability of the delayed optimal control is also fully characterized. For future research, the infinite-horizon stochastic LQ problem with input delay should be investigated.





\section*{Appendix}

\subsection*{A. Proof of Theorem \ref{Theorem-necessary}}\label{appendix-A}         

{
{$(i)\Rightarrow (ii)$}. Let $u^{t,x,*}$ be an optimal control of Problem (LQ) for the initial pair $(t,x)$. Then, we now prove that (\ref{Theorem-necessary-equality}) is satisfied with property (\ref{Z-X}). The following deduction is a variant of that in \cite{Zhang-huanshui2015}.

 Firstly, let us begin with the cast $k=N-1$. Noting $Z^{t,x,*}_{N}=GX_N^{t,x,*}$, we have
} %
\begin{eqnarray*}
\mathbb{E}_{N-1-d}Z^{t,x,*}_{N}
=G{A}_{N-1}\mathbb{E}_{N-1-d}X_{N-1}^{t,x,*}+G{B}_{N-1}u_{N-1}^{t,x,*},
\end{eqnarray*}
and
\begin{eqnarray*}
\mathbb{E}_{N-1-d}(Z^{t,x,*}_{N}w_{N-1})=G{C}_{N-1}\mathbb{E}_{N-1-d}X_{N-1}^{t,x,*}+G{D}_{N-1}u_{N-1}^{t,x,*}.
\end{eqnarray*}
Hence, (\ref{stationary-condition}) for $k=N-1$ reads as
\begin{eqnarray*}
&&0=R_{N-1} u^{t,x,*}_{N-1}+B_{N-1}^T\mathbb{E}_{{N-1}-d}Z^{t,x,*}_{{N}}+D_{N-1}^T\mathbb{E}_{{N-1}-d}(Z^{t,x,*}_{N}w_{N-1})\\[1mm]
&&\hphantom{0}=W_{N-1}u^{t,x,*}_{N-1}+H_{N-1}\mathbb{E}_{N-1-d}X^{t,x,*}_{N-1}.
\end{eqnarray*}
As there exists a $u^{t,x,*}$ satisfies (\ref{stationary-condition}), from Lemma \ref{Lemma-matrix-equation}  we know that (\ref{Theorem-necessary-equality}) holds for $k=N-1$, and that $u^{t,x,*}_{N-1}$ can be selected as
\begin{eqnarray*}
u^{t,x,*}_{N-1}=-W^\dagger_{N-1}H_{N-1}\mathbb{E}_{N-1-d}X^{t,x,*}_{N-1}.
\end{eqnarray*}
Furthermore,
\begin{eqnarray*}
&&Z^{t,x,*}_{N-1}=\big{(}Q_{N-1}+A_{N-1}^TGA_{N-1}+C_{N-1}^TGC_{N-1}\big{)} X^{t,x,*}_{N-1}-H_{N-1}^TW^\dagger_{N-1}H_{N-1}\mathbb{E}_{N-1-d}X^{t,x,*}_{N-1}\\[1mm]
&&\hphantom{Z^{t,x,*}_{N-1}}= P_{N-1}^{(0)} X^{t,x,*}_{N-1}+P^{(d)}_{N-1} \mathbb{E}_{N-1-d}X^{t,x,*}_{N-1}.
\end{eqnarray*}
In view of $P^{(i)}_{N-1}=0, i=1,\cdots, d-1$, we have (\ref{Z-X}) for $k=N-1$.

Secondly,  assume that for $k\in\{t+d,...,N-1\}$
\begin{eqnarray}\label{Theorem-necessary-proof-1-1}
&&H_\ell\mathbb{E}_{\ell-d}X_\ell^{t,x,*}\in \mbox{Ran}(W_\ell),~~\ell\in \mathbb{T}_{k+1}=\{k+1,...,N-1\},\\[1.5mm]
&&\label{Theorem-necessary-proof-1-2}
u_\ell^{t,x,*}=-W_\ell^\dagger H_\ell\mathbb{E}_{\ell-d}X^{t,x,*}_\ell,~~\ell\in \mathbb{T}_{k+1},
\end{eqnarray}
and
\begin{eqnarray*}
Z^{t,x,*}_{\ell+1}=P^{(0)}_{\ell+1}X_{\ell+1}^{t,x,*}+P^{(1)}_{\ell+1}\mathbb{E}_{\ell}X^{t,x,*}_{\ell+1}+\cdots+P^{(d)}_{\ell+1}\mathbb{E}_{{\ell+1}-d}X_{\ell+1}^{t,x,*},~~\ell\in \mathbb{T}_{k+1}.
\end{eqnarray*}
Now we verify that these are also true for the case  $k$. In fact, notice that
\begin{eqnarray}\label{Theorem-necessary-proof-1}
\mathbb{E}_{k-d}Z_{k+1}^{t,x,*}=\sum_{i=0}^{d}P^{(i)}_{k+1} \big{(}A_k\mathbb{E}_{k-d}X^{t,x,*}_k+B_ku_k^{t,x,*} \big{)},
\end{eqnarray}
and
\begin{eqnarray}\label{Theorem-necessary-proof-2}
\mathbb{E}_{k-d}(Z_{k+1}^{t,x,*}w_k)=P^{(0)}_{k+1}\big{(}C_k\mathbb{E}_{k-d}X^{t,x,*}_k+D_ku_k^{t,x,*} \big{)}.
\end{eqnarray}
Then, (\ref{stationary-condition}) reads as
\begin{eqnarray}\label{Theorem-necessary-proof-3}
&&0=R_k u^{t,x,*}_k+B_k^T\mathbb{E}_{k-d}Z^{t,x,*}_{k+1}+D_k^T\mathbb{E}_{k-d}(Z^{t,x,*}_{k+1}w_k)\nonumber \\
&&\hphantom{0}=W_ku_k^{t,x,*}+H_k\mathbb{E}_{k-d}X^{t,x,*}_k.
\end{eqnarray}
%
%
This implies by Lemma \ref{Lemma-matrix-equation} that (\ref{Theorem-necessary-equality}) holds for $k$ and that $u^{t,x,*}_k$ can be selected as
\begin{eqnarray}\label{Theorem-necessary-proof-4}
u_k^{t,x,*}=-W_k^\dagger H_k\mathbb{E}_{k-d}X^{t,x,*}_k.
\end{eqnarray}
Furthermore,
\begin{eqnarray}\label{Theorem-necessary-proof-5}
&&Z_k^{t,x,*}=\Big{[}Q_k+A_k^T\big{(}P^{(0)}_{k+1}+P^{(1)}_{k+1}\big{)}A_k+C_k^TP^{(0)}_{k+1}C_k \Big{]}X^{t,x,*}_k+\sum_{i=2}^{d} A_k^TP^{(i)}_{k+1}A_k\mathbb{E}_{k+1-i}X_k^{t,x,*}\nonumber \\
&&\hphantom{Z_k^{t,x,*}=}-H_k^TW_k^\dagger H_k\mathbb{E}_{k-d}X^{t,x,*}_k\nonumber \\[1mm]
&&\hphantom{Z_k^{t,x,*}}= P^{(0)}_kX_k^{t,x,*}+P^{(1)}_k\mathbb{E}_{k-1}X^{t,x,*}_{k}+\cdots+P^{(d)}_k\mathbb{E}_{k-d}X^{t,x,*}_k.
\end{eqnarray}

Let us further investigate a special case $k=t+d$ of (\ref{Theorem-necessary-proof-5})
\begin{eqnarray}\label{Theorem-necessary-proof-7}
Z_{t+d}^{t,x,*}=P^{(0)}_{t+d}X_{t+d}^{t,x,*}+P^{(1)}_{t+d}\mathbb{E}_{{t+d}-1}X^{t,x,*}_{{t+d}}+\cdots+P^{(d)}_{t+d}\mathbb{E}_{{t}}X^{t,x,*}_{t+d}.
\end{eqnarray}
Then, from a derivation similar to (\ref{Theorem-necessary-proof-1})-(\ref{Theorem-necessary-proof-4}), we have that (\ref{Theorem-necessary-equality}) holds for $k=t+d-1$ and
\begin{eqnarray*}
u_{t+d-1}^{t,x,*}=-W_{t+d-1}^\dagger H_{t+d-1}\mathbb{E}_{t}X^{t,x,*}_{t+d-1}.
\end{eqnarray*}
Therefore,
\begin{eqnarray}\label{Theorem-necessary-proof-8}
&&Z_{t+d-1}^{t,x,*}=\Big{[}Q_{t+d-1}+A_{t+d-1}^T\big{(}P^{(0)}_{{t+d}}+P^{(1)}_{{t+d}}\big{)}A_{t+d-1}+C_{t+d-1}^TP^{(0)}_{{t+d}}C_{t+d-1} \Big{]}X^{t,x,*}_{t+d-1}\nonumber \\
&&\hphantom{Z_{t+d-1}^{t,x,*}=}+\sum_{i=2}^{d-1} A_{t+d-1}^TP^{(i)}_{{t+d}}A_{t+d-1}\mathbb{E}_{{t+d}-i}X_{t+d-1}^{t,x,*}\nonumber \\
&&\hphantom{Z_{t+d-1}^{t,x,*}=}+\big{[}A_{t+d-1}^TP^{(d)}_{{t+d}}A_{t+d-1}-H_{t+d-1}^TW_{t+d-1}^\dagger H_{t+d-1}\big{]}\mathbb{E}_{t}X^{t,x,*}_{t+d-1}\nonumber \\[1mm]
&&\hphantom{Z_{t+d-1}^{t,x,*}}= P^{(0)}_{t+d-1}X_{t+d-1}^{t,x,*}+P^{(1)}_{t+d-1}\mathbb{E}_{{t+d-2}}X^{t,x,*}_{{t+d-1}}+\cdots+P^{(d-1)}_{t+d-1}\mathbb{E}_{t}X^{t,x,*}_{t+d-1}.
\end{eqnarray}
Note that the form of (\ref{Theorem-necessary-proof-8}) is different from (\ref{Theorem-necessary-proof-5}) and (\ref{Theorem-necessary-proof-7}).  Therefore, we further need deductions.

Assume that for $k\in \{t,...,t+d-2\}$ we have (\ref{Theorem-necessary-proof-1-1}), (\ref{Theorem-necessary-proof-1-2}) and
\begin{eqnarray*}
Z_{k+1}^{t,x,*}=P^{(0)}_{k+1}X_{{k+1}}^{t,x,*}+P^{(1)}_{k+1}\mathbb{E}_{k}X^{t,x,*}_{k+1}+\cdots+P^{(k+1-t)}_{{k+1}}\mathbb{E}_{{t}}X^{t,x,*}_{k+1}.
\end{eqnarray*}
%
%
%
%
%
%
%
%
Similar to (\ref{Theorem-necessary-proof-1})-(\ref{Theorem-necessary-proof-4}), we have that (\ref{Theorem-necessary-equality}) holds for $k$ and  $u^{t,x,*}_k$ can be selected as
\begin{eqnarray*}\label{Theorem-necessary-proof-6}
u_k^{t,x,*}=-W_k^\dagger H_k\mathbb{E}_{k-d}X^{t,x,*}_k.
\end{eqnarray*}
Furthermore,
\begin{eqnarray*}
&&Z_k^{t,x,*}=\Big{[}Q_k+A_k^T\big{(}P^{(0)}_{k+1}+P^{(1)}_{k+1}\big{)}A_k+C_k^TP^{(0)}C_k \Big{]}X_k^{t,x,*}+\sum_{i=2}^{k-t}A_k^TP^{(i)}_{k+1}A_k\mathbb{E}_{k+1-i}X_k^{t,x,*}\\[1mm]
&&\hpm{Z_k^{t,x,*}=}+\big{[}A_k^TP^{(k+1-t)}_{k+1}A_k-H_k^TW_k^\dagger H_k\big{]}\mathbb{E}_tX_k^{t,x,*}\\[1mm]
&&\hpm{Z_k^{t,x,*}}=P_k^{(0)}X^{t,x,*}_k+P_k^{(1)}\mathbb{E}_{k-1}X^{t,x,*}_k+\cdots+P^{(k-t)}_{k}\mathbb{E}_tX_k^{t,x,*}.
\end{eqnarray*}
By induction, we can prove (\ref{Theorem-necessary-equality}), (\ref{Theorem-necessary-u}) and (\ref{Z-X}).

{$(ii)\Rightarrow(i)$}. By Lemma \ref{Lemma-matrix-equation} and reversing the proof of {$(i)\Rightarrow(ii)$},  we can achieve the result. \endpf

\subsection*{B. Proof of Lemma \ref{add-subtr-1}}\label{appendix-B}

By adding to and subtracting
\begin{eqnarray*}
&&\sum_{k=t}^{N-1}\mathbb{E}\Big{\{}\sum_{i=0}^{d}(\mathbb{E}_{k+1-i}X^0_{k+1})^TP^{(i)}_{k+1}\mathbb{E}_{k+1-i}X^0_{k+1}-\sum_{i=0}^{d}(\mathbb{E}_{k-i}X^0_{k})^TP^{(i)}_{k}\mathbb{E}_{k-i}X^0_{k}\Big{\}}\\
&&+\sum_{k=t}^{t+d-1}\mathbb{E}\Big{\{}\sum_{i=0}^{k+1-t}(\mathbb{E}_{k+1-i}X^0_{k+1})^TP^{(i)}_{k+1}\mathbb{E}_{k+1-i}X^0_{k+1}-\sum_{i=0}^{k-t}(\mathbb{E}_{k-i}X^0_{k})^TP^{(i)}_{k}\mathbb{E}_{k-i}X^0_{k}\Big{\}}
\end{eqnarray*}
from $J(t,0;u)$, we have (noting $X^0_t=0$)
\begin{eqnarray}\label{add-subtr}
&&\hspace{-2em}J(t,0;{u}) =\sum_{k=t+d}^{N-1}\mathbb{E}\Big{\{}(X_k^0)^T\big{[}Q_k+A_k^T\big{(}P^{(0)}_{k+1}+P^{(1)}_{k+1}\big{)}A_k^T+C_k^TP^{(0)}_{k+1}C_k -P_{k}^{(0)}\big{]}X_k^0 \nonumber  \\
&&\hspace{-2em}\hphantom{J(t,0;{u}) =}+\sum_{i=1}^{d-1}(\mathbb{E}_{k-i}X^0_k)^T\big{[}A_k^TP^{(i+1)}_{k+1}A_k-P_k^{(i)} \big{]}\mathbb{E}_{k-i}X^0_k-(\mathbb{E}_{k-d}X^0_k)^TP^{(d)}_k\mathbb{E}_{k-d}X^0_k\nonumber  \\
&&\hspace{-2em}\hphantom{J(t,0;{u}) =}+2(H_k\mathbb{E}_{k-d}X^0_k)^Tu_k+u_k^TW_ku_k \Big{\}}\nonumber  \\
&&\hspace{-2em}\hphantom{J(t,0;{u}) =}+\sum_{k=t+2}^{t+d-1}\mathbb{E}\Big{\{}(X_k^0)^T\big{[}Q_k+A_k^T\big{(}P^{(0)}_{k+1}+P^{(1)}_{k+1}\big{)}A_k^T+C_k^TP^{(0)}_{k+1}C_k -P_{k}^{(0)}\big{]}X^0_k \nonumber  \\
&&\hspace{-2em}\hphantom{J(t,0;{u}) =}+\sum_{i=1}^{k-t-1}(\mathbb{E}_{k-i}X^0_k)^T\big{[}A_k^TP^{(i+1)}_{k+1}A_k-P_k^{(i)} \big{]}\mathbb{E}_{k-i}X^0_k\nonumber\\
&&\hspace{-2em}\hphantom{J(t,0;{u}) =}+(\mathbb{E}_{t}X^0_k)^T\big{(}A_k^TP^{(k+1-t)}_{k+1}A_k-P^{(k-t)}_k\big{)}\mathbb{E}_{t}X^0_k+2(H_k\mathbb{E}_{t}X^0_k)^Tu_k+u_k^TW_ku_k \Big{\}}\nonumber  \\[1mm]
&&\hspace{-2em}\hphantom{J(t,0;{u}) =}+\mathbb{E}\Big{\{}(X_{t+1}^0)^T\big{[}Q_{t+1}+A_{t+1}^T\big{(}P^{(0)}_{{t+2}}+P^{(1)}_{{t+2}}\big{)}A_{t+1}^T+C_{t+1}^TP^{(0)}_{{t+2}}C_{t+1} -P_{{t+1}}^{(0)}\big{]}X^0_{t+1} \nonumber  \\[1mm]
&&\hspace{-2em}\hphantom{J(t,0;{u}) =}+(\mathbb{E}_{t}X^0_{t+1})^T\big{(}A_{t+1}^TP^{(2)}_{t+2}A_{t+1}-P^{(1)}_{t+1}\big{)}\mathbb{E}_{t}X^0_{t+1}+2(H_{t+1}\mathbb{E}_{t}X^0_{t+1})^Tu_{t+1}\nonumber  \\[1mm]
&&\hspace{-2em}\hphantom{J(t,0;{u}) =}+u_{t+1}^TW_{t+1}u_{t+1} \Big{\}}+\mathbb{E}\Big{\{}(X_{t}^0)^T\big{[}Q_{t}+A_{t}^T\big{(}P^{(0)}_{{t+1}}+P^{(1)}_{{t+1}}\big{)}A_{t}^T+C_{t}^TP^{(0)}_{{t+1}}C_{t} -P_{{t}}^{(0)}\big{]}X^0_{t} \nonumber  \\[1mm]
%
%
%
&&\hspace{-2em}\hphantom{J(t,0;{u}) =}+2(H_{t}X^0_{t})^Tu_{t}+u_{t}^TW_{t}u_{t} \Big{\}}\nonumber  \\[1mm]
&&\hspace{-2em}\hpm{J(t,0;{u}) }=\sum_{k=t}^{N-1}\mathbb{E}\Big{\{}(\mathbb{E}_{k-d}X^0_k)^TH_k^TW_k^{\dagger}H_k\mathbb{E}_{k-d}X^0_k+2(H_k\mathbb{E}_{k-d}X^0_k)^Tu_k+u_k^TW_ku_k \Big{\}}.
%
%
\end{eqnarray}
This completes the proof. \endpf

\subsection*{C. Proof of Theorem \ref{Theorem-all-ii}}\label{appendix-C}

$(i)\Rightarrow(ii)(iii)$. Consider Problem (LQ) for the initial pair $(N-1,\xi)$ with $\xi\in \mathbb{R}^n$. Letting $k=N-1$ in (\ref{Riccati-1}), similar to (\ref{add-subtr}) we have
\begin{eqnarray}\label{Theorem-Nec-suff-all-1}
&&J(N-1,\xi; u_{N-1})=u_{N-1}^TW_{N-1}u_{N-1}+2(H_{N-1}\xi)^Tu_k+\xi^T H_{N-1}^TW_{N-1}^\dagger H_{N-1}\xi\nonumber\\
&&\hphantom{J(N-1,\xi; u_{N-1})=} +\xi^T\Big{(}\sum_{i=0}^dP^{(i)}_{N-1}\Big{)}\xi>-\infty.
\end{eqnarray}
As (\ref{Theorem-Nec-suff-all-1}) holds for any $\xi\in \mathbb{R}^n$ and any $u_{N-1}\in \mathcal{U}_{ad}^{N-1}$, we must have
\begin{eqnarray*}
W_{N-1}\geq 0,~~\mbox{Ran}(H_{N-1})\subset \mbox{Ran}(W_{N-1})\mbox{ (i.e.,} W_{N-1} W_{N-1}^\dagger H_{N-1}=H_{N-1}).
\end{eqnarray*}
Otherwise, if $W_{N-1}$ has a negative eigenvalue, say $\mu$, then for an eigenvector $\eta$ of $\mu$
\begin{eqnarray}\label{Theorem-Nec-suff-all-2}
J(N-1,\xi;\lambda \eta) = \mu\lambda^2 |\eta|^2+2\lambda (H_{N-1}\xi)^T\eta+\xi^T\Big{(}\sum_{i=0}^dP^{(i)}_{N-1}\Big{)}\xi\rightarrow -\infty,\mbox{ as }\lambda \rightarrow \infty.
\end{eqnarray}
This contradicts the finiteness of Problem (LQ). Further, if $\mbox{Ran}(H_{N-1})\subset \mbox{Ran}(W_{N-1})$, then there exists a $\xi_0\in \mathbb{R}^n$ such that $H_{N-1}\xi_0 \in\mbox{Ker}(W_k)$. Hence,
\begin{eqnarray*}
J(N-1,\xi;-\lambda H_{N-1}\xi_0) = -2\lambda |H_{N-1}\xi_0|^2+\xi^T\Big{(}\sum_{i=0}^dP^{(i)}_{N-1}\Big{)}\xi\rightarrow -\infty,\mbox{ as }\lambda \rightarrow \infty.
\end{eqnarray*}
This also contradicts the finiteness of Problem (LQ). We therefore have (\ref{Theorem-Nec-suff-all-2}) and
\begin{eqnarray}\label{Theorem-Nec-suff-all-3}
J(N-1, \xi; -W_{N-1}^\dagger H_k \xi)=\inf_{u_{N-1}}J(N-1,\xi;u_{N-1})=\xi^T\Big{(}\sum_{i=0}^dP^{(i)}_{N-1}\Big{)}.
\end{eqnarray}

Assume that for $k\in \mathbb{T}_{t+d}=\{t+d,...,N-1\}$
\begin{eqnarray}\label{Theorem-Nec-suff-all-4}
W_{\ell}\geq 0,~~W_{\ell} W_{\ell}^\dagger H_{\ell}=H_{\ell},~~\ell\in \mathbb{T}_{k+1}=\{k+1,...,N-1\},
\end{eqnarray}
and
\begin{eqnarray}\label{Theorem-Nec-suff-all-5}
J(\ell, \xi; u^{\ell,\xi,*})=\inf_{u\in \mathcal{U}_{ad}^{\ell}}J(\ell,\xi; u),~~\ell\in \mathbb{T}_{k+1},~ \xi\in \mathbb{R}^n.
\end{eqnarray}
In (\ref{Theorem-Nec-suff-all-5}), $u^{\ell,\xi,*}$ is given by
\begin{eqnarray*}
u_p^{\ell,\xi,*}=
-W_p^\dagger H_p\mathbb{E}_{p-d}X^{k,\xi,*}_p,~~p\in \mathbb{T}_{\ell}=\{\ell,...,N-1\}
\end{eqnarray*}
with
\begin{eqnarray*}
\left\{
\begin{array}{l}
X^{\ell,\xi,*}_{p+1}=\big{(}A_{p}X^{\ell,\xi,*}_p-B_{p}W_p^\dagger H_p\mathbb{E}_{p-d}X^{\ell,\xi,*}_p\big{)}+\big{(}C_{p}X^{\ell,\xi,*}_p-D_{p}W_p^\dagger H_p\mathbb{E}_{p-d}X^{\ell,\xi,*}_p\big{)}w_p,\\[1mm]
X^{\ell,\xi,*}_\ell=\xi,~~p\in \mathbb{T}_{\ell}.
\end{array}
\right.
\end{eqnarray*}
We now prove that (\ref{Theorem-Nec-suff-all-4}) and (\ref{Theorem-Nec-suff-all-5}) also hold for the case $\ell=k$. In fact, similarly to (\ref{add-subtr}) we have
\begin{eqnarray}\label{Theorem-Nec-suff-all-6}
&&J(k,\xi;u)=\sum_{\ell=k+1}^{N-1}\mathbb{E}\Big{\{}(u_\ell+W^\dagger_\ell H_\ell \mathbb{E}_{\ell-d}X_\ell)^TW_\ell(u_\ell+W_\ell^{\dagger}H_\ell\mathbb{E}_{\ell-d}X_\ell) \Big{\}}\nonumber\\
&&\hphantom{J(k,\xi;u)=}+\mathbb{E}\Big{\{}\xi^TH_k^TW_k^{\dagger}H_k\xi+2(H_k\xi)^Tu_k+u_k^TW_k u_k \Big{\}}+\xi^T\Big{(}\sum_{i=0}^dP^{(i)}_{k}\Big{)}\xi\nonumber\\[1mm]
&&\hphantom{J(k,\xi;u)}>-\infty,
\end{eqnarray}
which holds for any $\xi\in \mathbb{R}^n$ and any $u\in \mathcal{U}_{ad}^k$. Let the elements $u_{k+1},...,u_{N-1}$ of $u$  take the following form
\begin{eqnarray*}
u_\ell=
-W_\ell^\dagger H_\ell\mathbb{E}_{\ell-d}X_\ell,~~\ell\in \mathbb{T}_{k+1},
\end{eqnarray*}
and denote such a $u$ by $\widehat{u}$ with its element $u_k$ being freely selected. Then, (\ref{Theorem-Nec-suff-all-6}) becomes to
\begin{eqnarray}\label{Theorem-Nec-suff-all-7}
J(k,\xi;\widehat{u})=\mathbb{E}\Big{\{}\xi^TH_k^TW_k^{\dagger}H_k\xi+2(H_k\xi)^Tu_k+u_k^TW_k u_k \Big{\}}+\xi^T\Big{(}\sum_{i=0}^dP^{(i)}_{k}\Big{)}\xi >-\infty.
\end{eqnarray}
By an analysis similar to that between (\ref{Theorem-Nec-suff-all-2}) and (\ref{Theorem-Nec-suff-all-3}), we have
\begin{eqnarray*}
W_{k}\geq 0,~~W_{k} W_{k}^\dagger H_{k}=H_{k},
\end{eqnarray*}
and
\begin{eqnarray*}
J(k, \xi; u^{k,\xi,*})=\inf_{u\in \mathcal{U}_{ad}^{\ell}}J(k,\xi; u)=\xi^T\Big{(}\sum_{i=0}^dP^{(i)}_{k}\Big{)}\xi.
\end{eqnarray*}
with $u^{k,\xi,*}$ being given by (\ref{Theorem-all-u}).
By induction, we can get (ii) and (iii).

$(ii)\Rightarrow (i)$. This is straightforward.

$(ii)\Rightarrow (iii)$. By the proof of Theorem \ref{Theorem-necessary}, we know
\begin{eqnarray}\label{Theorem-all-ii-1}
0=W_\ell u^{k,\xi,*}_\ell+H_\ell \mathbb{E}_{\ell-d}X^{k,\xi,*}_\ell,&~~\ell\in \mathbb{T}_{k}.
\end{eqnarray}
Note that (\ref{Theorem-all-ii-1}) holds for any initial pair $(k,\xi)$ and $\ell\in \mathbb{T}_k$. We must have %
\begin{eqnarray*}
W_kW_k^\dagger H_k=H_k,~~k\in \mathbb{T}_t.
\end{eqnarray*}
%
%
%
$W_k\geq 0, k\in \mathbb{T}_t$, is due to the convexity of $u\mapsto J(t,x;u)$.

$(iii)\Rightarrow(ii)$. This is straightforward from
\begin{eqnarray*}
J(k,\xi;u)=\sum_{\ell=k}^{N-1}\mathbb{E}\Big{\{}(u_\ell+W^\dagger_\ell H_\ell \mathbb{E}_{\ell-d}X_\ell)^TW_\ell(u_\ell+W_\ell^{\dagger}H_\ell\mathbb{E}_{\ell-d}X_\ell) \Big{\}}+\Pi_k(\xi),
\end{eqnarray*}
where
\begin{eqnarray*}
\Pi_k(\xi)=\left\{
\begin{array}{l}
\sum_{i=0}^{k-t}\xi^TP_k^{(i)}\xi,~~~k\in \{t,...,t+d-1\},\\
\sum_{i=0}^{d}\xi^TP_k^{(i)}\xi,~~~k\in \mathbb{T}_{t+d}.
\end{array}
\right.
\end{eqnarray*}

$(iii)\Rightarrow(iv)$. Let $(P^{(0)},...P^{(d)})$ be the solution of Riccati-like equation set (\ref{Riccati-1})-(\ref{Riccati-2}) with property $W_kW_k^\dagger H_k=H_k, W_k\geq 0, k\in \mathbb{T}_t$. From the extended Schur's lemma, we know $(P^{(0)},...P^{(d)})\in \mathcal{M}$. Hence, $\mathcal{M}$ is nonempty.

$(iv)\Rightarrow(i)$. Let $(P^{(0)},...P^{(d)})\in \mathcal{M}$. Then, from Lemma \ref{Lemma-add-subtr}, we know
\begin{eqnarray*}
J(k,\xi;u)\geq \Pi_k(\xi)>-\infty.
\end{eqnarray*}
Hence, Problem (LQ) is finite.   \endpf




\end{document}